\documentclass[12pt]{article}

\usepackage{latexsym,amsmath,amsthm,amssymb,fullpage}

\usepackage{color}

\usepackage{latexsym,amsmath,amsthm,amssymb,fullpage}

\usepackage[utf8]{inputenc}

\usepackage{amssymb}
\usepackage{epstopdf}
\usepackage{amsmath}
\usepackage{amsthm}
\usepackage{amsfonts}
\usepackage{epstopdf}

\newtheorem{theorem}{Theorem}

\newtheorem{corollary}[theorem]{Corollary}

\newtheorem{proposition}[theorem]{Proposition}
\theoremstyle{definition}

\theoremstyle{remark}

\newcommand\beq{\begin {equation}}
\newcommand\eeq{\end {equation}}
\newcommand\beqs{\begin {equation*}}
\newcommand\eeqs{\end {equation*}}

\newcommand\R{{\mathbb R}}

\newcommand\PP{{\mathbb P}}
\newcommand\E{{\mathbb E}}

\newcommand{\mathbbm}[1]{\text{\usefont{U}{bbm}{m}{n}#1}}

\begin {document}

$$   $$

\vskip 4mm
\begin{center}
\textbf{\LARGE Four Talagrand inequalities}
\vskip 2mm
\textbf{\LARGE under the same umbrella}
\vskip 10mm
\textit{\large Michel Ledoux}\\
\vskip 3mm
{\large University of Toulouse, France}\\
\end{center}
\vskip 7mm

\vskip 4mm

\begin {abstract}
This note reviews the studies of the last decades emphasizing a common principle
based on entropy, logarithmic Sobolev inequality and
hypercontractivity, behind four most celebrated inequalities by M.~Talagrand:
the convex distance inequality, the $\mathrm {L}^1$--$\mathrm {L}^2$ variance inequality,
the quadratic transportation cost inequality, and the inequality on the supremum
of empirical processes.
\end {abstract}

\vskip 7mm

\section {Four Talagrand inequalities} \label {sec.1}

This section outlines the statements of the four Talagrand inequalities 
considered in this note, in the original notation of the author. 
These inequalities were all established (published) between 1991 and 1996,
the first three in relatively short articles.

\medskip

\noindent \textbf {Talagrand's convex distance inequality \cite [Theorem 1.1 and (1.3)] {T91},
\cite [Theorem~4.1.1]{T95}.}
Let $(\Omega_i , \mu_i )$, $i=1, \ldots , n$, be arbitrary
probability spaces and provide their product with the product probability $P$.
Given a set $ A \subset \Omega$ and $x \in \Omega$, set
$$
U_A (x) \, = \, \big \{ {(s_i)}_{1\leq i \leq n} \in \{0,1\}^n ; \exists \, y \in A ;
		s_i = 0 \Rightarrow x_i = y_i \big \}.
$$
Denote by $V_A(x)$ the convex hull of $U_A(x)$ considered as a subset of $\R^n$.
The set $V_A(x)$ contains $0$ if and only if $x$ belongs to $A$. Denote then by
$d_A(x)$ the Euclidean distance of $0$ to $V_A(x)$ (in \cite {T95}, the notation
$f_c(A,x)$ is used instead, the letter $c$ refereeing to ``convexity").

For any (measurable) $ A \subset \Omega$,
\beq \label {eq.convexhull}
	\int_\Omega e^{\frac 14  d_A^2} dP \, \leq \, \frac {1}{P(A)} \, .
\eeq

\medskip

\noindent \textbf {Talagrand's $\mathrm {L}^1$--$\mathrm {L}^2$ variance inequality
\cite [Theorem 1.5]{T94}.}
Let $ \{0,1\}^n$ be the discrete cube equipped with the product probability measure $\mu_\mathsf{p} $
giving weight $\mathsf{p}$ to $1$ and $1-\mathsf{p}$ to $0$, $0 < \mathsf{p} < 1$.
For $1 \leq r \leq \infty$, denote by ${\| \cdot \| }_r$ the norm in $\mathrm {L}^r(\mu_\mathsf{p} )$.
For every $i = 1, \ldots, n$ and every $ x = (x_1, \ldots, x_n) \in  \{0, 1\}^n$, let
$U_i(x)$ be the point of $ \{0,1\}^n$ obtained from $x$ by replacing $x_i$ by $1-x_i$ and leaving
the other coordinates unchanged. If $f$ is a function on $ \{0,1\}^n$, set
$ \Delta_i f(x) = (1-\mathsf{p}) (f(x) - f(U_i(x)))$ if $x_i = 1$ and
$ \Delta_i f(x) = \mathsf{p} (f(x) - f(U_i(x)))$ if $x_i = 0$.

There is a numerical constant $K>0$ such that for any function $f : \{0,1\}^n \to \R$
such that $\int_{\{0,1\}^n} f d\mu_\mathsf{p} = 0$,
\beq \label {eq.l1l2}
	{\| f\|}_2^2 \, \leq \,  K \, \log \Big ( \frac {2}{\mathsf{p}(1-\mathsf{p})} \Big) 
      \sum_{i=1}^n \frac {{\| \Delta_i f\| }_2^2} 
          {\log \big ( e \,\frac {{\| \Delta_i f\| }_2}{{\| \Delta_i f\| }_1} \big ) } \, .
\eeq

\medskip

\medskip

\noindent \textbf {Talagrand's quadratic transportation cost inequality
\cite [Theorem 1.1] {T96a}.}
For two probability measures $\nu$ and $\mu$ on the Borel sets of $\R^n$, define the transportation
cost $T_w (\mu, \nu) $ as the infimum of
$$
\int_{\R^n \times \R^n} w(x,y) d\pi (x,y)
$$
over all the probability measures $\pi $ on $\R^n \times \R^n$
such that $\mu$ is the first marginal of $\pi$ and $\nu$ is the second, and where
$w(x,y) = \sum_{i=1}^n (x_i - y_i)^2$.
Let $d\gamma (x) = e^{-|x|^2/2} \frac {dx}{(2\pi)^{n/2}}$
be the standard Gaussian measure on the Borel sets of $\R^n$.

For every probability measure
$\mu$ absolutely continuous with respect to $\gamma$ with $f = \frac {d\mu}{d\gamma}$,
\beq \label {eq.transport}
T_w (\mu ,\gamma ) \, \leq \,  2 \int_{\R^n} f \log f \, d\gamma 
		\, = \, 2 \int_{\R^n} \log f \, d\mu.
\eeq

\medskip

\noindent \textbf {Talagrand's inequality on the supremum of empirical processes 
\cite [Theorem 1.4] {T96b}.}
Let $X_1, \ldots, X_n$ be independent random variables with values in a
measurable space $\Omega$, and let $\mathcal {F}$ be a countable
family of measurable functions on $\Omega$. Consider the random variable 
$$
Z \, = \, \sup_{f \in  {\cal F}} \sum _{i=1}^n f(X_i).
$$
Set $U = \sup_{f \in \mathcal {F}} {\| f\|}_\infty$ and
$V = \E \big ( \sup_{f \in  {\cal F}} \sum _{i=1}^n f(X_i)^2 \big)$. 

For each $t >0$, 
\beq \label {eq.supremum}
  \PP \big (  \big |Z -  \E (Z)\big | \geq  t \big )
  \, \leq  \,  K \exp \bigg ( -  \frac {t}{KU} \log \bigg ( 1 + \frac{tU}{V} \bigg ) \bigg )
\eeq
where $K>0$ is a numerical constant.

\section {Introduction} \label {sec.2}

M.~Talagrand's mathematical achievements have deeply influenced the scientific
developments over the last decades. His work embraces numerous mathematical
fields, including measure theory, functional analysis, Banach space geometry, stochastic processes,
Boolean analysis, isoperimetric and concentration inequalities, optimal transportation, statistical physics etc.

It is in particular extremely impressive to hear about \textit {four} of his
inequalities, all called \textit {Talagrand's inequality}, in rather different areas.
The precise statements of these inequalities have been presented in the first section.

\noindent -- 
The famous convex distance inequality for product measures, put forward in his celebrated Publication
of the Institut des Hautes \'Etudes Scientifiques \cite {T95}, has very fruitful applications
in numerous areas of probability theory, statistics or optimization, and combinatorial and
discrete mathematics.

\noindent -- The $\mathrm {L}^1$--$\mathrm {L}^2$ variance inequality on the discrete cube
first provided an alternate, sharper, approach to the Kahn-Kalai-Linial theorem on influences
in Boolean analysis, and became a central tool in theoretical
computer science. As put forward by I.~Benjamini, G.~Kalai and O.~Schramm, 
it is besides, at this point, the only generic argument towards
sub-diffusive and super-concentration phenomena which are ubiquitous to many models of the current research
in probability theory (percolation, random matrices, spin glasses etc.).

\noindent -- Talagrand's quadratic transportation cost has been one
founding stone in the interaction between partial differential equations, probability theory and geometry, 
as first emphasized by F.~Otto and C.~Villani in connection
with the logarithmic Sobolev inequality, leading to the new field
of functional inequalities, curvature lower bounds and 
analysis on metric measure spaces by J.~Lott, C.~Villani, K.~T.~Sturm, L.~Ambrosio, N.~Gigli, G.~Savaré.

\noindent -- The methods developed for concentration inequalities for product measures
led to the fundamental Talagrand
inequality for the supremum of empirical processes, a major and essential tool in
modern infinite-dimensional statistics.

\noindent The Talagrand inequalities have both shaped the mathematics of the last decades and
stand today as essential and common tools of the current research.

\medskip

Apart perhaps the $\mathrm {L}^1$--$\mathrm {L}^2$ variance inequality, the Talagrand
inequalities were mostly motivated by aspects of the concentration of measure phenomenon
and its applications to Banach space geometry and probability in Banach spaces
\cite {MS86,P86,LT91,G98,L01,Sc03,AGM15}.
``The idea of concentration of measure, which was discovered by V.~Milman,
is arguably one of the great ideas of analysis in our times'' (M.~Talagrand \cite {T96c}).
The concentration of measure phenomenon
has indeed become today a major tool in various areas of mathematics such as
asymptotic geometric analysis, probability theory,
statistical mechanics, mathematical statistics and learning theory, random matrix theory or
quantum information theory, stochastic dynamics, randomized algorithms, complexity...
Of isoperimetric origin and flavour, it is
suited to the investigation of models involving an infinite number of variables, and
emphasizes that,
quoting M.~Talagrand again \cite {T96c}, ``a random variable that depends (in a ``smooth'' way)
on the influence of many independent random variables (but not too much on any of them)
is essentially constant''.
It is indeed a main feature of the four Talagrand inequalities that they are dimension free
(constants do not depend on the size of the samples, and the statements
extend to infinite-dimensional systems).
These inequalities are in particular inspired by, and 
provide deep and powerful extensions of, the model concentration inequality for Gaussian measures
expressing that for any Lipschitz function $F : \R^n \to \R$ with
Lipschitz semi-norm ${\| F\|}_{\mathrm {Lip}}$,
\beq \label {eq.gaussianconcentration}
\gamma \big ( F \geq \textstyle {\int_{\R^n} F d\gamma} + r \big)
		\, \leq \, e^{-r^2/2 {\|F\|}_{\mathrm {Lip}}^2}
\eeq
for all $r \geq 0$ (where $d\gamma (x) = e^{-|x|^2/2} \frac {dx}{(2\pi)^{n/2}}$,
the standard Gaussian measure on the Borel sets of $\R^n$).

The Talagrand inequalities look quite different, and have been established by him with
different tools and methods (although induction on the dimension towards
dimension-free bounds may be detected as a common background). It is
the purpose of this note to show, based on the mathematical developments of the last decades,
that all these four inequalities may be seen as consequences of a common principle, namely entropy,
logarithmic Sobolev inequality and hypercontractivity. 

This observation started with the papers \cite {DS84} and \cite {AMS94}
which revived an unpublished letter by I.~Herbst to L.~Gross in 1975
deducing exponential integrability
from the logarithmic Sobolev inequality, just after the fundamental discovery by L.~Gross of
the latter \cite {G75}.
The contribution \cite {AMS94} by S.~Aida, T.~Masuda and I.~Shigekawa has been very
influential in this regard. The relevance of this observation towards
the Gaussian concentration inequality \eqref {eq.gaussianconcentration} was pointed out next in \cite {L95}.
Based on this principle, the contribution \cite {L96b} provided
an alternate approach to the inequality on supremum of empirical processes \eqref {eq.supremum},
later developed
and extended in several steps and contributions summarized in the monograph \cite {BLM13}
by S.~Boucheron, G.~Lugosi and P.~Massart. A few years later, these authors also covered the
convex distance inequality \eqref {eq.convexhull} by this method. 
(The monograph \cite {BLM13} actually covers in depth various aspects emphasized in this note.)
Already present in \cite {T94} and earlier \cite {KKL88}, the hypercontractivity character
of the $\mathrm {L}^1$--$\mathrm {L}^2$ variance inequality \eqref {eq.l1l2} is clarified in
\cite {BKS03}, and led these authors to a proof of sub-linearity of percolation times.
It is also mentioned in \cite {BKS03}
that Talagrand's variance inequality \eqref {eq.l1l2} may be used to recover the associated
concentration inequality for percolation time first obtained in \cite {T95}, linking even more
the four inequalities. In the early 2000,
F.~Otto and C.~Villani \cite {OV00} deduced the Talagrand quadratic transportation cost
inequality from the logarithmic Sobolev inequality of L.~Gross,
bringing to light the deep link between these two families of inequalities.

The purpose of this note is thus to provide a common, self-contained, 
approach to the four Talagrand inequalities based
on entropy, logarithmic Sobolev inequality and hypercontractivity, with for each of them a simple
and direct proof. Sub-additivity properties of entropy and 
logarithmic Sobolev inequality will encapsulate the dimension-free concentration phenomenon,
and the Talagrand inequalities.

The next paragraph,
Section~\ref {sec.3}, presents the necessary (elementary) material on the entropy method to this task.
It is restricted to the purpose of the subsequent proofs, but wide extensions and applications
of the approach have been developed in the literature. The proofs of the
Talagrand inequalities are addressed in the following four sections. It should be mentioned that, with
respect to the original Talagrand's formulations, 
the notation and terminology will somewhat be adapted in order to follow the common trends in this regard.
This is also motivated by the framework emphasized in Section~\ref {sec.3}.
No attempt is made towards sharp numerical constants (besides the one in
\eqref {eq.transport}), favouring the simplicity of the arguments before technical improvements
(with appropriate references for the latter).
Each section will be followed by some comments on historical aspects
and further developments, with a few (and only a few) general pointers to the relevant literature.

This note should have been written several years ago. Hopefully it is still of some interest.

\section {Entropy, logarithmic Sobolev inequality, hypercontractivity} \label {sec.3}

This section describes the basic material towards the entropy approach to the four Talagrand inequalities.
As already mentioned, the presentation of the various tools is limited to this specific task, but, as emphasized in the literature, the power and generality of the principle may be extended and applied far outside
this given purpose. This material is rather elementary, and most of the statements may be found
in standard monographs or lecture notes (see the notes and references at the end of the paragraph)
involved with these objects.

$(\Omega, \Sigma, \mu)$ will denote a generic probability space. For each $1 \leq p \leq \infty$,
${\| \cdot \|}_p$ is the norm of the Lebesgue space $\mathrm {L}^p (\mu)$.
The various integrability conditions appearing below will be automatically satisfied in all the
subsequent illustrations dealing mostly with bounded functions.

\medskip

\noindent \textbf {Entropy.} For a measurable
function $f : (\Omega, \Sigma, \mu) \to \R$, non-negative in $\mathrm {L}^1(\mu)$, set
$$
\mathrm {Ent}_\mu (f) \, = \, \int_\Omega f \log f \,  d\mu 
		- \int_\Omega f \, d\mu \, \log \bigg ( \int_\Omega f \, d\mu  \bigg)
		\, \in \, [0 , + \infty]
$$
(since $u \in [0,\infty) \mapsto u\log u$ is convex bounded from below, with the convention
$0 \log 0 = 0$).
Observe that $\mathrm {Ent}_\mu (f)$ is homogeneous of order $1$.
If $d\nu = f d\mu$ for a density $f$,
$\mathrm {Ent}_\mu (f) = \mathrm {H}(\nu \, | \,\mu)$, the relative entropy of $\nu$
with respect to $\mu$. 

Entropy has a lot of common with variance, 
$$
\mathrm {Var}_\mu (f) \, = \, \int_\Omega f^2  d\mu  - \bigg (\int_\Omega f \, d\mu \bigg)^2 
$$
(for a function $f $ in $\mathrm {L}^2(\mu)$), in particular the following duality and variational
representations.

The duality formula for entropy expresses that
\beq \label {eq.entropyduality} 
\mathrm {Ent}_\mu (f)  \, = \, \sup \bigg \{ \int_\Omega fg \, d\mu \, ;
		\int_\Omega e^g d\mu  \leq 1 \bigg   \} 
\eeq
(where $g$ may be assumed to be bounded from above and below).
Indeed, assume by homogeneity that $\int_\Omega \! fd\mu  =1$. By Young's inequality
$ uv \leq   u\log u - u + e^v$, $u \geq 0$, $v \in  \R$,
so that, for $\int_\Omega e^g d\mu  \leq 1$,
$$ 
\int_\Omega fg \, d\mu  \,\leq \,   \int_\Omega f\log f \, d\mu  -1 +  \int_\Omega e^g d\mu 
                                  \leq  \int_\Omega f\log f \, d\mu  \, = \, \mathrm {Ent}_\mu (f).
$$
For the converse, set $f_N = \min (\max (f, \frac 1N), N)$, $N \geq 1$, and choose
$g = \log \big (\frac {f_N}{\int_\Omega f_N d\mu}\big)$. The claim follows as $N \to \infty$.

For the further purposes, note that this duality formula justifies the well-known
entropic inequality, for any $f \geq 0$ with $\int_\Omega f d\mu = 1$ and $g$ measurable
such that $fg$ is integrable,
\beq \label {eq.entropicinequality}
\int_\Omega fg \, d\mu \, \leq \, 
		\mathrm {Ent}_\mu (f) + \log \bigg (\int_\Omega e^g d\mu \bigg) 
		\, = \, 	\int_\Omega f\log f \, d\mu + \log \bigg (\int_\Omega e^g d\mu \bigg).
\eeq

The variational formula on the other hand states that
\beq \label {eq.entropyvariational} 
\mathrm {Ent}_\mu (f)  
		\, = \, \inf_{c>0} \int_\Omega \big [ f(\log f - \log c) - (f-c) \big] d\mu.
\eeq
Indeed, the infimum of $ c  \mapsto  c - (\log c + 1) \int_\Omega f d\mu $
is attained at $c = \int_\Omega f d\mu$ giving thus rise to $\mathrm {Ent}_\mu (f)$.

\medskip

\noindent \textbf {Tensorization of entropy.}
A fundamental feature of entropy (and of variance) is its product or sub-additivity
property, main source of the approach to the Talagrand inequalities developed here.
Let $(\Omega _i, \Sigma  _i, \mu _i)$,
$i=1, \ldots ,n$, be probability spaces, and denote by $P$ the product probability measure
$\mu _1 \otimes \cdots \otimes \mu _n$ on the product space
$X = \Omega _1 \times \cdots \times \Omega _n$ equipped with the product
$\sigma $-field. A point $x$ in $X$ is denoted
$x = (x_1, \ldots ,x_n)$, $x_i \in \Omega _i$, $i=1, \ldots ,n$.
Given $f$ on the product space, write furthermore $f_i$, $i=1, \ldots ,n$,
for the function on $\Omega _i$ defined by
$$f_i(x_i) = f(x_1, \ldots , x_{i-1}, x_i, x_{i+1}, \ldots ,x_n),$$
with $x_1, \ldots , x_{i-1}, x_{i+1}, \ldots ,x_n$ fixed.

\begin {proposition} \label {prop.tensorization}
For every non-negative function $f$ on the product space $X$ in $\mathrm {L}^1(P)$,
$$
\mathrm {Ent}_P (f) 	\, \leq  \, \sum _{i=1}^n \int_X \mathrm {Ent}_{\mu _i} (f_i) dP .
$$
\end {proposition}

It may be pointed out that on the right-hand side, integration in $dP$ is actually,
for each $i=1, \ldots, n$, over the remaining
coordinates $(x_1, \ldots, x_{i-1}, x_{i+1}, \ldots, x_n)$.

\begin {proof}
Use the duality formula \eqref {eq.entropyduality}. Given $g$
such that $\int_X e^g dP\leq  1$, set, for every $i = 1 \ldots ,n$,
$$ 
g^i (x_i, \ldots , x_n) 
	\, = \,  \log \bigg ( \frac{ \int_X e^{g(x_1, \ldots ,x_n)}d\mu _1(x_1) \cdots 
                        d\mu _{i-1}(x_{i-1}) } { \int_X e^{g(x_1, \ldots ,x_n)}d\mu _1(x_1) \cdots 
                        d\mu _{i}(x_{i}) } \bigg )
$$
(well-defined for $\mu_i \otimes \cdots \otimes \mu_n$-almost every $(x_i, \ldots , x_n)$).
Then $ g \leq  \sum _{i=1}^n g^i$ and $\int_{\Omega_i} e^{(g^i)_i} d\mu _i = 1$
for every $i = 1 \ldots ,n$. Therefore,
\beqs \begin {split}
	\int_X fg \, dP 
           & \, \leq \,  \sum _{i=1}^n \int_X fg^i \,dP \\
	&\, = \, \sum _{i=1}^n \int_X \bigg (\int_{\Omega_i}  f_i(g^i)_i d\mu _i \bigg )dP \\
           & \, \leq \, \sum _{i=1}^n \int_X \mathrm{Ent}_{\mu _i} (f_i) dP \\
\end {split} \eeqs
which is the result.
\end {proof}

Proposition~\ref {prop.tensorization} is presented in \cite {L96b}
(in a more general form due to S.~Bobkov, the proof presented here being due
to S.~Kwapien), and deduced from Han's inequality in \cite {BLM13}.
It has a classical analogue for the variance, established in the
same (even simpler) way, known as the Efron-Stein inequality (cf.~\cite {ES81, S86, RT86, S97,BLM13})
expressing that for every function $f$ on the product space $X$ in $\mathrm {L}^2(P)$,
\beq \label {eq.efronstein}
 \mathrm {Var}_P (f)  \, \leq \,  \sum _{i=1}^n \int_X \mathrm {Var}_{\mu _i} (f_i) dP.
\eeq
Actually \eqref {eq.efronstein} may be deduced from Proposition~\ref {prop.tensorization}
applied to $f = 1 + \varepsilon g$ with $\varepsilon \to 0$.

The tensorization Proposition~\ref {prop.tensorization} admits related formulations which
will be of fundamental use in the applications to the Talagrand convex
distance inequality and the inequality on the supremum of empirical processes.

First, by Jensen's inequality, for every $i=1,\ldots, n$ (and $x_1, \ldots, x_{i-1}, x_{i+1}, \ldots, x_n$
fixed),
$$
\mathrm{Ent}_{\mu _i} (f_i) \, \leq \, \frac 12
		\int_{\Omega_i} \int_{\Omega_i} \big [f_i(x_i) - f_i(y_i) \big]
				\big [\log f_i(x_i) - \log f_i(y_i) \big] d\mu_i(x_i) d\mu_i(y_i)
$$
(an analogue of the duplication formula for the variance). Therefore,
at a first (mild) level, Proposition~\ref {prop.tensorization} yields by symmetry
\beq \begin {split} \label {eq.tensorization1}
& \mathrm {Ent}_P (f) 	\\
& \, \leq  \, 
	\sum _{i=1}^n \int_X   \int \! \int_{\{f_i(x_i) \geq f_i(y_i) \}} 
		 \big [f_i(x_i) - f_i(y_i) \big]
				\big [\log f_i(x_i) - \log f_i(y_i) \big] d\mu_i(x_i) d\mu_i(y_i)  dP(x) \\
\end {split} \eeq

Proposition~\ref {prop.tensorization} may also be combined with the variational representation
\eqref {eq.entropyvariational} in the form
\beq  \label {eq.tensorization2}
\mathrm {Ent}_P (f) 	
	\, \leq  \, \sum _{i=1}^n \int_X 
	\bigg (\inf_{c_i > 0} \int_{\Omega_i} 
	\big [ f_i (\log f_i - \log c_i) - (f_i  -c_i) \big ] \, d\mu_i \bigg) dP . 
\eeq
Of course, for each $i = 1, \ldots, n$, $c_i > 0$ in the infimum may be chosen to depend on the variables
$x_1, \ldots, x_{i-1}, x_{i+1}, \ldots, x_n$.

\medskip

\noindent \textbf {Logarithmic Sobolev inequality on the two-point space.}
A logarithmic Sobolev inequality bounds the entropy of a function $f$ by an energy, or Dirichlet
form. It is the analogue of the classical Poincaré inequality for the variance. It has been discovered
by L.~Gross \cite {G75} in 1975 together with its equivalence with hypercontractivity
(although earlier versions may detected in the literature).

The most basic logarithmic Sobolev inequality takes place on the two-point space
$\{-1, +1 \}$ with the Bernoulli probability
measure $\mu_\mathsf{p} (\{+1\}) = \mathsf{p}$,  $\mu_\mathsf{p} (\{-1\}) = 1- \mathsf{p} = \mathsf{q}$, 
$0 < \mathsf{p} < 1$, and states that for any $f : \{-1, +1 \} \to \R$,
\beq \label {eq.logsob2point}
\mathrm {Ent}_{\mu_\mathsf{p}} (f^2) 
\, \leq  \, \frac 1\rho \, \mathsf{p} \mathsf{q}  \big [ f(+1) - f(-1) \big ]^2
\eeq
where $ \rho = \frac {\mathsf{p}-\mathsf{q}}{ \log \mathsf{p} - \log \mathsf{q}} $. 
The symmetric case $\mathsf{p = \frac 12}$ is the simplest and most classical one, but as will be seen,
general values of $\mathsf{p}$ are handled below similarly.

When $\mathsf{p = \frac 12}$, $\rho = \frac 12$ so that
the constant on the right-hand side of \eqref {eq.logsob2point} is $\frac 12$.
A proof in this case runs as follows.
Setting $ f(+1) = \alpha$ and $f(-1) = \beta$, the inequality amounts to
$$
\frac {\Phi (\alpha^2) + \Phi (\beta^2)}{2} - \Phi \Big ( \frac {\alpha^2 + \beta^2}{2} \Big)
		\, \leq \,  \frac 12 \, (\alpha - \beta)^2
$$
where $\Phi (u) = u\log u$, $u \geq 0$.
That is, if $ r = \frac 12 (\alpha^2 + \beta^2) $ and $ s = \frac 12 (\alpha^2 - \beta^2) $
$$
\Phi (r+s) + \Phi (r-s) - 2 \Phi (r) \, \leq \,  (\alpha - \beta)^2.
$$
But the left-hand side is
$$
\int_0^s \big [ \Phi '' (r+v) + \Phi ''(r-v ) \big] (b-v) dv
$$
and since the function $\Phi'' = \frac 1u$ is convex on the given domain
$$
\Phi '' (r+v) + \Phi ''(r-v ) \, \leq \, 2 \Phi''(r) \, = \, \frac 2r \, .
$$
It remains to observe that $\frac {s^2}{r} \leq (\alpha - \beta)^2$.
A direct proof for general $\mathsf{p}$ may be found in \cite {ABC+00}.

\medskip

\noindent \textbf {Hypercontractivity on the two-point space.}
The logarithmic Sobolev inequality \eqref {eq.logsob2point} 
may be translated equivalently into the famous hypercontractivity property.
At this stage, only the symmetric case $\mathsf{p=\frac 12}$ is described for simplicity.
Any function  $ f : \{ -1, +1\} \to \R$ may be represented as $f(x) = a + bx$,
$x \in \{ -1, +1\}$, $a, b \in \R$. For any $t\geq 0$, define the new function
$ P_t f(x) = a +  e^{-t} b x$, $x \in \{ -1, +1\}$. The family
$ {(P_t)}_{t \geq 0}$ defines a semigroup of contractions on $\mathrm {L}^p(\mu_{\mathsf{\frac 12}})$
for any $1 \leq p \leq \infty$. It turns out that \eqref {eq.logsob2point} 
allows for the strengthening
\beq \label {eq.hypercontractivity}
{\| P_t f \|}_q \, \leq \, { \|  f \|}_p 
\eeq
whenever $1 < p <q < \infty$, $ e^{2t} \geq   \frac {q-1}{p-1} $.
The latter \eqref {eq.hypercontractivity} may be translated again as a two-point inequality
$$
\Big ( \frac 12 \, |a + e^{-t} b|^q + \frac 12 \, |a - e^{-t} b|^q  \Big)^{1/q}
		\, \leq \, \Big ( \frac 12 \, |a +  b|^p + \frac 12 \, |a -  b|^p  \Big)^{1/p} .
$$
This inequality was established directly by A.~Bonami \cite {B71} and W.~Beckner \cite {B75}
but L.~Gross \cite {G75} observed that it is actually equivalent to \eqref {eq.logsob2point} 
(and as such easier to establish). This connection is developed next, as well as the case
$\mathsf{p \not= \frac 12}$, after 
the setting is extended to the product model $X = \{-1,+1\}^n$, $n \geq 1$.

\medskip

\noindent \textbf {Logarithmic Sobolev inequality and hypercontractivity on the discrete cube.}
On the two-point space $\{-1,+1\}$ equipped with the Bernoulli measure $\mu_\mathsf{p}$,
$0 < \mathsf{p} < 1$, consider the (Markov)
operator $\mathrm {L} f = \int_{\{-1,+1\}} f d\mu_\mathsf{p} - f$. Note that
\beq \label {eq.dirichlet1}
\int_{\{-1,+1\}} f(- \mathrm {L}f) d\mu_\mathsf{p} 
	\, = \,  \int_{\{-1,+1\}} (\mathrm {L}f)^2 d\mu_\mathsf{p}
		\, = \, \mathrm {Var}_{\mu_\mathsf{p}} (f) 
		\, = \, \mathsf{p} \mathsf{q} \big [ f(+1) - f(-1) \big]^2.
\eeq

On the product space $X = \{-1,+1\}^n$ with the product measure 
$\mu^n_\mathsf{p}$, consider the product
operator $\mathrm {L} = \sum_{i=1}^n \mathrm {L}_i$ where
$\mathrm {L}_i$ is acting on the $i$-th coordinate of a function $f : X \to \R$
as $\mathrm {L}_i f = \int_{\{-1,+1\}} f_i(x_i) d\mu_\mathsf{p}(x_i) - f$, $i = 1, \ldots, n$
(with the notation of the tensorization paragraph, that is 
$f_i(x_i) = f(x_1, \ldots, x_{i-1}, x_i, x_{i+1}, \ldots, x_n)$ with
$x_1, \ldots, x_{i-1}, x_{i+1}, \ldots, x_n$ fixed).

The Dirichlet form
$$
\mathcal {E} (f,g) \, = \, \int_X f (-\mathrm {L} g) d\mu_\mathsf{p}^n
$$
for two functions $f, g : X \to \R$ admits various representations of the form 
\beq \begin {split} \label {eq.dirichlet2}
& \mathcal {E}  (f,g)  \\
	& \, = \, \sum_{i=1}^n \int_X \mathrm {L}_i f \, \mathrm {L}_i g \, d\mu_\mathsf{p}^n \\
	& \, = \, \sum_{i=1}^n \int_X \mathrm {Cor}_{\mu_\mathsf{p}} (f_i , g_i) d\mu_\mathsf{p}^n \\
	& \, = \, \frac 12 \sum_{i=1}^n \int_X 
		\int_{\{-1, +1\}} \! \int_{\{-1, +1\}} \big [ f_i(x_i) - f_i(y_i) \big ] \big [ g_i(x_i) - g_i(y_i) \big ]
				d\mu_\mathsf{p}(x_i) d\mu_\mathsf{p}(y_i) d\mu_\mathsf{p}^n (x).
\end {split} \eeq
In particular, from the Efron-Stein inequality \eqref {eq.efronstein}, the Poincaré inequality
\beq \label {eq.poincare}
\mathrm {Var}_{ \mu_\mathsf{p}^n} (f)  \, \leq \, \mathcal {E} (f,f)
\eeq
holds true for any $f : X \to \R$.

The operator $\mathrm {L}$ generates the Markov semigroup ${(P_t)}_{t\geq 0}$ defined by
\beq \label {eq.semigroupcube}
P_t \, = \, \sum_{k = 0}^\infty \frac {t^k}{k !} \,  \mathrm {L}^k
\eeq
($\mathrm {L}^0 = \mathrm {Id}$)
in the sense that $\frac {d}{dt} P_t f = \mathrm {L}P_tf =  P_t \, \mathrm {L}f$.
It is symmetric and invariant with respect to $\mu_\mathsf{p}^n$, that is
for functions $f, g : X \to \R$, $\int_X f P_t g \, d\mu_\mathsf {p}^n 
= \int_X g P_t f \, d\mu_\mathsf{p}^n$. It is immediately checked that for
$\mu_{\mathsf{\frac 12}}$ on $\{-1,+1\}$, $P_t f(x) = a + e^{-t}bx$
for a function $f = a + bx$, $x \in \{-1,+1\}$.

It is a consequence of the tensorization property, and the representations
\eqref {eq.dirichlet1} and \eqref {eq.dirichlet2} of the associated Dirichlet forms,
that the logarithmic Sobolev inequality
\eqref {eq.logsob2point} extends to functions $f$ on the product space
$X = \{-1, +1 \}^n$ equipped with the product measure $\mu_\mathsf{p}^n$ in the form
\beq \label {eq.logsobcube}
 	  \mathrm {Ent}_{\mu_\mathsf{p}^n} (f^2)
        \, \leq  \,  \frac 1\rho \, \mathcal {E} (f,f)
\eeq
for any $f : X \to \R$, where it is recalled that
$ \rho = \frac {\mathsf{p} - \mathsf{q}}{\log \mathsf{p} - \log \mathsf{q}}$
($=\frac 12$ if $\mathsf{p = \frac 12}$).
In the same way, the hypercontractivity inequality \eqref {eq.hypercontractivity} extends
to functions $f$ on $ X = \{-1, +1 \}^n$, expressing that whenever
$1 < p < q < \infty$ and $e^{4\rho  t} \geq \frac {q-1}{p-1}$,
\beq \label {eq.hypercontractivitycube}
{\| P_t f \|}_q \, \leq \, { \|  f \|}_p \, .
\eeq

While the tensorization of hypercontractivity may be achieved independently (cf.~e.g. \cite {B75,B94}), 
it is fruitful
to deduce it from the logarithmic Sobolev inequality as developed by L.~Gross \cite {G75}.
Given $1 < p < q < \infty$, the key idea of L.~Gross is to differentiate in time the quantity
${\| P_t f\|}_{q(t)}$, $q(t) = (p-1) e^{4\rho  t} - 1$,
$t \geq 0$. Since $|P_t f| \leq P_t(|f|)$, it is enough to deal with a non-negative
function $f$ (actually not identically zero). It holds that
$$
\frac {d}{dt} \int_X (P_t f)^{q(t)} d\mu_\mathsf{p}^n
	\, = \, q'(t) \int_X (P_t f)^{q(t)} \log P_t f \, d\mu_\mathsf{p}^n
			+ q(t) \int_X (P_t f)^{q(t)-1} \, \mathrm {L} P_t f \, d\mu_\mathsf{p}^n
$$
and hence
\beqs \begin {split}
{\| P_t f\|}^{q(t)-1}_{q(t)} & \,  \frac {d}{dt} {\| P_t f\|}_{q(t)} \\
	& \, = \,  - \frac {q'(t)}{q(t)^2} \int_X (P_t f)^{q(t)} d\mu_\mathsf{p}^n \, 
		 \log \int_X (P_t f)^{q(t)} d\mu_\mathsf{p}^n
			+ \frac {1}{q(t)} \frac {d}{dt} \int_X (P_t f)^{q(t)} d\mu_\mathsf{p}^n \\
	& \, = \, \frac {q' }{q^2}\, \mathrm {Ent}_{\mu_\mathsf{p}^n} \big ((P_t f)^q \big)
			+ \int_X (P_t f)^{q-1} \, \mathrm {L} P_t f \, d\mu_\mathsf{p}^n	 \\
\end {split} \eeqs
where the short-hand notation $q = q(t)$, $q' = q'(t)$, is used in the last line.
Assume now that
\beq \label {eq.convexity}
- \int_X (P_t f)^{q-1} \, \mathrm {L} P_t f \, d\mu_\mathsf{p}^n \, = \,
\mathcal {E} \big ( (P_t f)^{q-1}, P_tf \big) 
     \, \geq \,  \frac {4(q-1)}{q^2} \, \mathcal {E} \big ( (P_t f)^{q/2}, (P_t f)^{q/2} \big)
\eeq
so that 
$$
 q^2 \, {\| P_t f\|}^q_q \,  \frac {d}{dt} {\| P_t f\|}_{q} 
	 \, \leq \, q'\, \mathrm {Ent}_{\mu_\mathsf{p}^n} \big ((P_t f)^q \big)
		- 4(q-1)\, \mathcal {E} \big ( (P_t f)^{q/2}, (P_t f)^{q/2} \big).
$$
Applying the logarithmic Sobolev inequality \eqref {eq.logsobcube} to $(P_t f)^{q(t)/2}$
then indicates that for the choice of $q = q(t) = (p-1)e^{4\rho  t} -1$, 
the right-hand side of the latter inequality is negative. 
Therefore the map $ t \mapsto {\| P_t f\|}_{q(t)}$ is decreasing, which amounts to the hypercontractivity
inequality \eqref {eq.hypercontractivitycube}. The proof shows in the same way that
hypercontractivity is actually equivalent to the logarithmic Sobolev inequality.

It remains nevertheless to establish \eqref  {eq.convexity} which follows from a convexity argument.
Namely, for all $u > v \geq 0$,
\beqs \begin {split}
\bigg ( \frac {u^{q/2} - v^{q/2}}{u-v} \bigg)^2
	& \, = \, \bigg (\frac {q}{2(u-v)} \int_v^u s^{\frac q2-1} ds \bigg)^2\\
	& \, \leq \, \frac {q^2}{4(u-v)} \int_v^u s^{q-2} ds \\
	& \, = \, \frac {q^2}{4(q-1)} \, \frac {u^{q-1} - v^{q-1}}{u-v} \, .
\end {split} \eeqs
Hence, for any $u,v \in \R$,
$$
\big (u^{q-1} - v^{q-1} \big )(u-v) \, \geq \,  \frac {4(q-1)}{q^2} \, \big (u^{q/2} - v^{q/2}\big)^2.
$$
Recalling \eqref {eq.dirichlet2} then concludes to \eqref  {eq.convexity}.

\medskip

\noindent \textbf {Logarithmic Sobolev inequality for the Gaussian measure.}
If $d\gamma (x) = e^{-|x|^2/2} \frac {dx}{(2\pi)^n}$ is the standard Gaussian
(product) measure on the Borel sets of $\R^n$ ($|x|$ being the Euclidean length of $x \in \R^n$),
the logarithmic Sobolev inequality for $\gamma$ states that
\beq \label {eq.logsobgaussian}
\mathrm {Ent}_\gamma (f^2) \, \leq  \, 2 \int_{\R^n} |\nabla f|^2  d\gamma 
\eeq
for every smooth (for example locally Lipschitz in $\mathrm {L}^2(\gamma)$)
function $f : \R^n \to \R$. It might be worthwhile mentioning
that by the tensorization Proposition~\ref {prop.tensorization}, it is enough to know
the inequality in dimension one.

A proof of this inequality, also due to L.~Gross \cite {G75}, 
may be obtained from the logarithmic Sobolev inequality 
on the discrete cube via the central limit theorem. Roughly speaking,
if $\varphi : \R \to \R$ is smooth with compact support, apply
\eqref {eq.logsobcube} with $\mathsf {p=\frac12}$ to
$ f(x_1, \ldots, x_n) = \varphi \big (\frac {1}{\sqrt n} \sum_{i=1}^n x_i \big)$,
$ (x_1, \ldots, x_n) \in \{-1,+1\}^n$. If 
$ {\| \varphi' \|}_\infty +{\| \varphi'' \|}_\infty < \infty$, by a Taylor
expansion
$$
		f_i(x_i) - f_i(y_i) \, = \, \frac {x_i -y_i}{\sqrt n} \,
		\varphi' \bigg (\frac {1}{\sqrt n} \sum_{j\not= i} x_j \bigg) + O\Big ( \frac 1n \Big)
$$
uniformly in $x_i, y_i$, $i = 1, \ldots, n$. Hence,
$$
	 \sum_{i=1}^n \big [ f_i(x_i) - f_i(y_i) \big]^2
	 \, =\, \frac {1}{n} \sum_{i=1}^n (x_i -y_i)^2 \, 
		\varphi' \bigg (\frac {1}{\sqrt n} \sum_{j\not= i} x_j \bigg)^2 
		+  O\Big ( \frac {1}{\sqrt n} \Big)
$$
and from \eqref {eq.dirichlet2},
$$
\mathcal {E} (f,f) \, \leq \,
		 2 \int_X \varphi' \bigg (\frac {1}{\sqrt n} \sum_{j=2}^n x_j \bigg)^2 
		 d\mu_{\mathsf {\frac 12}}^n +  O\Big ( \frac {1}{\sqrt n} \Big).
$$
By the central limit theorem, $\limsup_{n \to \infty} \mathcal {E} (f,f) \leq 2\int_\R {\varphi'}^2 d\gamma$
and similarly $\lim_{n \to \infty} \mathrm {Ent}_{\mu_{\mathsf{\frac 12}}^n}(f^2)
= \mathrm {Ent}_\gamma (\varphi^2)$.

In addition to this proof, there are at least 15 different further
proofs of the logarithmic Sobolev inequality for the Gaussian measure.
The following introduces the analytic (semigroup) proof put forward by
D.~Bakry and M.~\'Emery \cite {BE85,B94} 
(which has been recognized as the simplest one by L.~Gross in 2010). 
It is presented with the Ornstein-Uhlenbeck semigroup (with invariant measure~$\gamma$); a
similar argument may be developed with the standard heat (Brownian) kernel/semigroup
(cf.~\cite {BGL14}).

Denote by ${(P_t)}_{t \geq 0}$ the so-called Ornstein-Uhlenbeck semigroup defined by
\beq \label {eq.ou}
P_t f(x) \, = \, \int_{\R^n} f \big ( e^{-t} x + \sqrt {1 - e^{-2t}} \, y \big) d\gamma (y),
		\quad t \geq 0, \, \,  x \in \R^n
\eeq
(for any $f : \R^n \to \R$ in $\mathrm {L}^1(\gamma)$ for example).
The Ornstein-Uhlenbeck semigroup is
invariant and symmetric with respect to the standard Gaussian measure $\gamma$,
with infinitesimal generator $\mathrm {L} f (x) = \Delta f (x) - \langle x , \nabla f (x) \rangle$
for which the integration by parts formula
$$
\int_{\R^n} f (-\mathrm {L}g) d\gamma \, = \, \int_{\R^n} \langle \nabla f , \nabla g \rangle \, d\gamma
$$
holds true for every smooth functions $f, g : \R^n \to \R$.

Let $f: \R^n \to \R$ be a (measurable) function such that
$\varepsilon \leq f \leq \frac 1\varepsilon$ for some $\varepsilon >0$ and
$\int_{\R^n} f d\gamma = 1$. For every $t>0$, $P_tf$ is then a $C^\infty$ function,
$\varepsilon \leq P_t f \leq \frac 1\varepsilon$,
and, as $t \to \infty$, $P_t f \to \int_{\R^n} f d\gamma = 1$ at every point. Therefore
$$
\int_{\R^n} f \log f \, d\gamma 
	\, = \, - \int_0^\infty \bigg (\frac {d}{dt} \int_{\R^n} P_t f \log P_t f \, d\gamma \bigg) dt.
$$
By the chain rule
\beqs \begin {split}
\frac {d}{dt} \int_{\R^n} P_t f \log P_t f \, d\gamma
	& \, = \, \int_{\R^n} \mathrm {L} P_t f \log P_t f \, d\gamma
			+  \int_{\R^n} \mathrm {L} P_t f \, d\gamma \\
	& \, = - \, \int_{\R^n} \frac {|\nabla P_t f|^2}{P_t f} \, d\gamma \\
\end {split} \eeqs
where integration by parts is used as well as the fact that $\int_{\R^n} \mathrm {L} g d\gamma = 0$.
Next, at any point, $|\nabla P_t f| \leq e^{-t} P_t (|\nabla f|)$ as is clear from the integral representation
\eqref {eq.ou} of $P_t$, and by the Cauchy-Schwarz inequality along the same representation,
$$
 \frac {|\nabla P_t f|^2}{P_t f} 
	\, \leq \, e^{-2t} \, \frac {P_t(|\nabla f|)^2}{P_t f}
	\, \leq \, e^{-2t} P_t \bigg (\frac {|\nabla f|^2} {f} \bigg) .
$$
Hence, by invariance of $P_t$ with respect to $\gamma$,
$$
  - \frac {d}{dt} \int_{\R^n} P_t f \log P_t f \, d\gamma
		\, \leq \, e^{-2t} \, \int_{\R^n} \frac {|\nabla  f|^2}{ f} \, d\gamma
$$
from which it follows by integration that
\beq \label {eq.logsobgaussian2}
 \int_{\R^n} f \log f \, d\gamma  
        \, \leq  \, \frac {1}{2} \int_{\R^n} \frac {|\nabla f|^2}{f} \, d\gamma.
\eeq
This is a classical alternate formulation of the logarithmic Sobolev inequality with,
on the left-hand side the relative entropy of $f d\gamma$ with respect to $\gamma$,
and on the right-hand side the so-called Fisher information of $f$.

Let now $f : \R^n \to \R$ be differentiable with gradient in $\mathrm {L}^2(\gamma)$.
Apply \eqref {eq.logsobgaussian2} to, for example, 
$$
\frac {(P_t f_N)^2 + \varepsilon}{\int_{\R^n} (P_t f_N)^2 d\gamma + \varepsilon}
$$
where $\varepsilon , t >0$ and $f_N = \min ( \max (f, -N), N))$. Letting successively
$\varepsilon \to 0$ and $t \to 0$ yields that
$$
\mathrm {Ent}_\gamma (f_N^2) \, \leq \, 2 \int_{\R^n} |\nabla f|^2 d\gamma
$$
uniformly in $N\geq 1$. By means of the entropic inequality \eqref {eq.entropicinequality}, for
every $\theta \in (0,1)$,
$$
(1-\theta) \int_{\R^n}f_N^2 d\gamma
	\, \leq \, \int_{\R^n} f_N^2 \mathbbm{1}_B \, d\gamma
	\, \leq \, 2 \int_{\R^n} |\nabla f|^2 d\gamma
					+ \int_{\R^n}f_N^2 d\gamma \, \log \big (1 + (e-1)\gamma (B) \big)
$$
where $B = \big \{ f_N^2 \geq \theta \int_{\R^n}f_N^2 d\gamma \big \}$, from which it easily
follows that $\sup_N \int_{\R^n}f_N^2 d\gamma < \infty$, and therefore by monotone convergence
that $\int_{\R^n}f^2 d\gamma < \infty$. Then also $\int_{\R^n} f^2 \log f^2 d\gamma < \infty$
which, altogether, concludes to  \eqref {eq.logsobgaussian}.

The constant $2$ is optimal in \eqref {eq.logsobgaussian} as can be seen from the choice
of the functions $f(x)= e^{\langle b,x \rangle}$, $x \in \R^n$, $b \in \R^n$, which achieve
equality.

On the basis of the logarithmic Sobolev inequality \eqref {eq.logsobgaussian} for the
Gaussian measure, the Gross argument may be developed similarly along the
Ornstein-Uhlenbeck semigroup to show that it is hypercontractive, with parameter $\rho = \frac 12$
in the notation of \eqref {eq.hypercontractivitycube}, a property going back to
E.~Nelson \cite {N66,N73} in quantum field theory. (It may be pointed out that in this
diffusive setting, the convexity argument \eqref {eq.convexity} is immediate by 
integration by parts and the chain rule.)

\medskip

\noindent \textbf {The Herbst argument.}
The main inspiration connecting a logarithmic Sobolev inequality to a concentration property
is the Herbst argument. It was originally observed by I.~Herbst towards exponential integrability
properties, with, in the notation below, the application of the Gross logarithmic Sobolev
inequality for the Gaussian measure
to $e^{\lambda F^2}$, $\lambda \in \R$. It is developed here on
$e^{\lambda F}$, $\lambda \in \R$, with in this form direct application to concentration
inequalities. This observation is at the root of the entropic proof of the Talagrand
convex distance inequality and on the supremum of empirical processes, as well
actually as the transportation cost inequality.
It is presented here on the Gaussian model (see the references for much more).

Let $F$ be a smooth bounded Lipschitz function on $\R^n$ with
Lipschitz semi-norm ${\| F\|}_{\mathrm {Lip}} \leq 1$.
Since $F$ is assumed to start with to be regular enough, it can be that $|\nabla F| \leq 1$ at every point.
The aim is to apply the logarithmic Sobolev inequality
\eqref {eq.logsobgaussian} to $f^2=  e^{\lambda F} $ for every $\lambda \in \R$.
With the notation $\Lambda (\lambda ) = \int_{\R^n} e^{\lambda F} d\gamma$,
$\lambda \in \R$, it holds that
$$
\mathrm {Ent}_\gamma (f^2)
    \, = \, \lambda\int_{\R^n} F e^{\lambda F} d\gamma 
    		-  \Lambda (\lambda) \log \Lambda(\lambda) 
		\, = \, \lambda \Lambda'(\lambda) - \Lambda (\lambda ) \log \Lambda (\lambda )
$$
while
$$
\int_{\R^n} |\nabla f|^2 d\gamma
		\, = \, \frac {\lambda^2}{4} \int_{\R^n} |\nabla F|^2 e^{\lambda F} d\gamma
		\, \leq \, \frac {\lambda^2}{4} \, \Lambda (\lambda).
$$
Hence from \eqref {eq.logsobgaussian},
$$ 
	 \lambda \Lambda '(\lambda )  -  \Lambda(\lambda)\log \Lambda(\lambda)
    \,    \leq  \, \frac {\lambda^2}{2} \, \Lambda(\lambda).
$$
If $H(\lambda ) = \frac 1\lambda  \log \Lambda (\lambda )$
(with $H(0) = \frac {\Lambda '(0)}{\Lambda(0)}= \int_{\R^n} F d\gamma $), $\lambda  \in \R$,
then $ H'(\lambda ) \leq \frac 12$ for every $\lambda $.  Therefore,
$  H(\lambda) - H(0)  \leq  \frac {\lambda}{2}$ if $\lambda \geq 0$
and $ \geq  \frac {\lambda}{2}$ if $\lambda \leq 0$
from which
\beq \label {eq.herbst}
\Lambda(\lambda) = \int_{\R^n} e^{\lambda F} d\gamma 
		\, \leq \, e^{\lambda \int_{\R^n}F d\gamma + \frac {\lambda^2}{2}}
\eeq
for every $\lambda \in \R$.

If $F$ is an arbitrary Lipschitz function with ${\| F\|}_{\mathrm {Lip}} \leq 1$,
apply the preceding for example to $P_t(F_N)$, $t>0$, $N \geq 1$,
where $P_t$ is the Ornstein-Uhlenbeck semigroup and
$F_N = \min ( \max (F, -N), N))$, and let then $t\to 0$ and $N \to \infty$ in \eqref {eq.herbst}. 

From \eqref {eq.herbst}, let $F : \R^n \to \R$ be Lipschitz with
Lipschitz semi-norm ${\| F\|}_{\mathrm {Lip}}$.
For every $r \geq 0$ and $\lambda \geq 0$, by Markov's inequality,
$$
\gamma \big (F \geq \textstyle {\int_{\R^n} F d\gamma} + r \big)
		 \, \leq \, e^{-\lambda (\int_{\R^n} F d\gamma + r)} \int_{\R^n} e^{\lambda F} d\gamma
			\, \leq \, e^{-\lambda r + \lambda^2/2 {\| F\|}_{\mathrm {Lip}}^2 } 
$$
where \eqref {eq.herbst} is used after homogeneity.
Optimizing in $\lambda$, the Gaussian concentration inequality \eqref {eq.gaussianconcentration}
$$
\gamma \big (F \geq \textstyle {\int_{\R^n} F d\gamma} + r \big)
		 \, \leq \, e^{- r^2/2 {\| F\|}_{\mathrm {Lip}}^2 } , \quad r \geq 0,
$$
is recovered in this way. Together with the same inequality for $-F$ and the union bound,
$$
\gamma \big ( | F - \textstyle {\int_{\R^n} F d\gamma}| \geq  r \big)
		 \, \leq \, 2 \, e^{- r^2/2 {\| F\|}_{\mathrm {Lip}}^2 }, \quad r \geq 0. 
$$
These inequalities describe the fundamental concentration property
of Gaussian measures (cf.~e.g. \cite {P86,LT91,Li95,L96a,B98,L01}).

\medskip

\textbf {Some notes and references.} As mentioned at the beginning, this section
puts forward the entropy method, and logarithmic Sobolev and hypercontractivity inequalities,
only in the context required by the proofs of the Talagrand inequalities
emphasized in this note. The framework and methodology may be
vastly extended and generalized to various settings and applications. Introductions
to logarithmic Sobolev inequalities are for example \cite {DS89, G93,R07,L99,ABC+00,GZ03,BGL14}...
The general monographs and courses from the non-exhaustive list \cite {DSC96,L01,V03,BLM13,OD14,GS15}...
contain further material and suitable pointers to the literature relevant to the topics of this section.

\section {The convex distance inequality} \label {sec.4}

In order to address the convex distance inequality via the entropy method
displayed in the previous section, it is of interest to first suitably translate
the (somewhat obscure) functional $d_A$ of \eqref {eq.convexhull}.

Recall the framework of a product probability measure
$P = \mu_1 \otimes \cdots \otimes \mu_n$
on the product space $X = \Omega_1 \times \cdots \times \Omega_n$.
A point $x$ in $X$ has coordinates $x= (x_1, \ldots, x_n)$.
There might be some measurability questions
in the forthcoming claims and proofs, but as emphasized in \cite {T95} these are unessential,
and one should ``treat all sets and functions as if they were measurable. This is certainly the case 
if one should assume that the $\Omega_i$'s are Polish and the $\mu_i$'s are Borel measures, and 
that one studies only compact sets, which is the only situation that occurs in applications".
If necessary, it is therefore even possible to assume all sets finite.

The weighted Hamming distance on $X$, with weight 
$a = (a_1, \ldots, a_n) \in  [0,\infty)^n$, is defined by 
$$
d_a(x,y) \, = \, \sum _{i=1}^n a_i \, \mathbbm {1}_{\{x_i \not= y_i\} } , \quad x, y \in X.
$$
For every non-empty (measurable) subset $A$ of $X $ and every $x \in X$, let then
$$
F_A(x) \, = \,  \sup _{|a|=1}d_a(x,A)
$$
where $|a| = \big (\sum_{i=1}^n a_i^2 \big)^{1/2}$.

The point is that the functional $F_A$ is actually equal to $d_A $ of \eqref {eq.convexhull}
\cite [Lemma 4.1.2] {T95}. Recall that, for $x \in X$,
$$
U_A (x) \, = \, \big \{ s = {(s_i)}_{1\leq i \leq n} \in \{0,1\}^n ; \exists \, y \in A ;
		s_i = 0 \Rightarrow x_i = y_i \big \}
$$
and that $d_A(x)$ is the Euclidean distance of $0$ to the convex hull $V_A(x)$ of $U_A(x)$
in $\R^n$. It is easily seen that, in this definition, $U_A(x)$ may be replaced
by the collection of the indicator functions $\mathbbm {1}_{\{x_i\neq y_i\} }$, $y \in  A$.
Now, if $d_A(x) <  r$ for some $r > 0$, there exists
$z$ in $V_A(x)$ with $|z| < r$. Let $a  \in  [0,\infty)^n$
with $|a|=1$. Then
$$
\inf _{y \in  V_A(x)} \langle a,y\rangle  \leq  \langle a,z\rangle  \, \leq \,  |z| \, < \,  r.
$$
Since
\beq \label {eq.technical}
\inf _{y \in  V_A(x)} \langle a,y\rangle  
	\, = \,  \inf _{s \in  U_A(x)} \langle a,s\rangle  \, = \,  d_a(x,A),
\eeq
it follows that $F_A(x) < r$. Hence $d_A(x) \geq F_A(x)$.
Conversely, assume that $d_A(x) >0$ (otherwise there will be nothing to prove) and
let $\delta >0$. Take $z \in V_A(x)$ such that $0 < |z|^2 \leq d_A(x)^2 + \delta$.
By convexity, for every $\theta \in (0,1)$ and every $y \in V_A(x)$,
$\theta y + (1-\theta ) z \in V_A(x)$ so that
$$ 
\big | z + \theta (y-z) \big |^2 
     \, = \,  \big | \theta y + (1-\theta ) z \big |^2 \, \geq \,  d_A(x)^2 \, \geq \, |z|^2 - \delta.
 $$
Therefore
$$
2 \theta \langle y-z,  z \rangle  + \theta^2 |y-z|^2 \, \geq \, - \delta.
$$
Setting $a = \frac {z}{|z|}$,
$$
 \langle a,y \rangle  \, \geq \, |z| - \frac {\theta |y-z|^2}{2|z|} - \frac { \delta}{2\theta |z|} \, \geq \,
 	|z| - \frac {2\theta n}{|z|} - \frac { \delta}{2\theta |z|}  \, .
$$
Now, by \eqref {eq.technical},
$$  
	F_A(x) \, \geq \, d_a(x,A) \, = \, \inf _{y \in V_A(x)}  \langle a ,  y \rangle 
	\, \geq \,  |z| - \frac {2\theta n}{|z|} - \frac { \delta}{2\theta |z|}  
	\, \geq \,  d_A(x) - \frac {2\theta n}{d_A(x)} - \frac { \delta}{2\theta d_A(x)}  \, .
$$
Since $\delta >0$ and $\theta \in (0,1)$ are arbitrary, it follows that $F_A(x) \geq d_A(x)$.

On the basis of this description,
the principle will be to apply the tensorization Proposition~\ref {prop.tensorization} 
to the functional $F_A$ together with the Herbst argument. For each $d_a(\cdot, A)$, or
any Lipschitz function with respect to the Hamming metric $d_a(x,y)$, a simple
Laplace transform argument together with iteration along the coordinates yields
Gaussian concentration bounds (cf.~\cite {MS86,T95,M98,L01}). 
The challenge is to achieve a similar goal uniformly over all $a$'s in the unit sphere of $\R^n$.

To ease the notation, set $F = F_A$ throughout the following steps
(assuming that ${P(F>0) >0}$ otherwise there is nothing to prove). Let $\varepsilon >0$.
For each $x = (x_1, \ldots, x_n) \in X$, 
there exists $a(x) = a = (a_1, \ldots, a_n) \in [0,\infty)^n$ with $|a|=1$ such that
$$
F(x) \, \leq \, d_a(x,A) + \varepsilon.
$$
For $1 \leq i \leq n$ and $y_i \in \Omega_i$,
set $y = (x_1, \ldots, x_{i-1}, y_i, x_{i+1}, \ldots, n)$. Then, with the notation of 
Proposition~\ref {prop.tensorization},
$$
F_i(x_i) - F_i(y_i) \, = \, F(x) - F(y) \, \leq \, d_a(x,A) - d_a(y,A) + \varepsilon.
$$
By the triangle inequality,
$$
F_i(x_i) - F_i(y_i) \, \leq \, d_a(x,y) + \varepsilon
		\, = \, a_i  \, \mathbbm {1}_{\{x_i \not= y_i\}} + \varepsilon \, \leq \, a_i + \varepsilon.
$$

Apply now Proposition~\ref {prop.tensorization} 
to $f = e^{\lambda F^2}$, $\lambda \geq 0$, in the form of \eqref {eq.tensorization1}.
Whenever $F_i(x_i) \geq F_i(y_i)$ ($\geq 0$) for $i=1, \ldots, n$, by the
mean value inequality and the preceding,
\beqs \begin {split}
\big [\lambda F_i(x_i)^2 -  \lambda F_i(y_i)^2 \big] 
	\big [ e^{\lambda F_i(x_i)^2} - e^{\lambda F_i(y_i)^2}\big] 
	 & \, \leq \, \lambda^2 \big [ F_i(x_i)^2 - F_i(y_i)^2 \big]^2  e^{\lambda F_i(x_i)^2} \\
	  & \, \leq \, 4 \lambda^2 \big (a_i(x) + \varepsilon \big)^2 F_i(x_i)^2  e^{\lambda F_i(x_i)^2} .\\
\end {split} \eeqs
As $\sum_{i=1}^n a_i(x)^2 = 1$, it therefore follows from \eqref {eq.tensorization1} that
$$
\mathrm {Ent}_P (e^{\lambda F^2}) 
		\, \leq \, 4 \lambda^2 \int_X 
		\big (1 + 2\varepsilon \sqrt {n}  + \varepsilon^2 n \big ) F^2 e^{\lambda F^2} dP
$$
for any $\lambda \geq 0$. As $\varepsilon \to 0$,
\beq \label {eq.herbsthull}
\mathrm {Ent}_P (e^{\lambda F^2}) 
		\, \leq \, 4 \lambda^2 \int_X F^2 e^{\lambda F^2} dP.
\eeq

Set $\Lambda (\lambda) = \int_X e^{\lambda F^2} dP$, $\lambda \in \R$, so that the preceding
expresses that, in the range $\lambda \geq 0$,
$$
  \lambda \Lambda '(\lambda) - \Lambda (\lambda) \log \Lambda (\lambda)
 		\, \leq \, 4\lambda^2 \Lambda '(\lambda).
$$
In the spirit of the Herbst argument, it remains to integrate such a differential inequality. 
If $K(\lambda ) = \log \Lambda (\lambda) $ ($ > 0$ when $\lambda >  0$),
it reads for $0 < \lambda < \frac 14$,
$$
\lambda (1 - 4\lambda) K'(\lambda) \, \leq \, K(\lambda),
$$
that is $(\log K)'(\lambda) \leq \frac {1}{\lambda (1 - 4\lambda)}$.
Hence, for every $0 < \eta \leq \lambda < \frac 14$,
$$
\log K(\lambda) - \log K(\eta)
		\, \leq \, \int_\eta^\lambda \frac {1}{u(1-4u)} \, du
		\, = \, \log \Big ( \frac {\lambda}{1 - 4\lambda} 
				\cdot \frac {1 -4\eta}{\eta} \Big).
$$
As $\eta \to 0$, $K(\eta) =  \log \int_X e^{\eta F^2} dP \sim \eta M_2$ where
$M_2 = \int_X F^2 dP$. Therefore, in this limit, the preceding inequality yields
$$
\log K(\lambda) - \log M_2 \, \leq \, \log \Big ( \frac {\lambda}{1 -4\lambda} \Big),
$$
that is
$$
 K(\lambda ) = \log \bigg (\int_X e^{\lambda F^2} dP\bigg)
		 \, \leq \, \frac {\lambda M_2}{1 -4\lambda} 
$$
for every $ 0 \leq \lambda < \frac 14$. For instance with $\lambda = \frac {1}{14}$ (see below
for the significance of this choice),
\beq \label {eq.convexhullintermediate1}
\int_X e^{\frac {1}{14} F^2} dP  \, \leq \, e^{\frac {1}{10} M_2}.
\eeq

As $\lambda \to 0$ in \eqref {eq.herbsthull}, $\mathrm {Var}_P(F^2) \leq 8M_2$.
Since $F = F_A = 0$ on $A$ it follows that $M_2 \leq \frac {8}{P(A)}$.
Hence, from \eqref {eq.convexhullintermediate1},
\beq \label {eq.convexhullintermediate2}
\int_X e^{\frac {1}{14} F^2} dP  \, \leq \,  e^{\frac {4}{5 P(A)}}.
\eeq
In particular, if $P(A) \geq \frac 12$, 
$\int_X e^{ \frac {1}{14} F_A^2} dP  \leq 5$.

These first conclusions are weaker than \eqref {eq.convexhull} in terms of numerical constants, but more
importantly in terms of the dependence on $P(A)$, especially for sets $A$ with
small probability. It is nevertheless already good enough for many of the significant applications
of Talagrand's convex distance inequality. It may be pointed out also that, towards
\eqref {eq.convexhullintermediate2}, the Herbst argument may be developed more
simply with $e^{\lambda F}$ rather than $e^{\lambda F^2}$. The presentation here is
motivated by homogeneity with the second part which is coming next.

To reach better dependence on $P(A)$ as in \eqref {eq.convexhull}, it is necessary
to develop the previous analysis but for negative values of $\lambda$.
Consider therefore \eqref {eq.tensorization1} now applied to 
$f = e^{-\lambda F^2}$, $\lambda \geq 0$.
To this task, it is worthwhile to observe that given $x = (x_1, \ldots, x_n)$ and $y_i \in \Omega_i$,
$i = 1, \ldots, n$,
\beq \label {eq.lipschitzsquare}
F_i(x_i)^2 - F_i(y_i)^2 \, = \, F(x)^2 - F(y)^2 \, \leq \, 1
\eeq
where as usual $F_i(y_i) = F (x_1, \ldots, x_{i-1}, y_i, x_{i+1}, \ldots, x_n) = F(y)$.
A proof of this claim may be given using the identity $F(x) = F_A(x) = d_A(x)$.
Indeed, 
$$
d_A(x)^2 
	 \, = \, \inf \big \{ |w|^2 ; \, w \in \mathrm {Conv} 
		\big( {(\mathbbm{1}_{\{x_j \not= z_j\} })}_{1 \leq j \leq n} ; z \in A\big) \big \} 
	 \, = \, \inf \bigg | \sum \theta_z s(z) \bigg|^2
$$
where the infimum is over all finite sums $\sum \theta_z s(z)$ with $\theta_z \geq 0$,
$ \sum \theta_z = 1$, $ s(z) = {(\mathbbm{1}_{\{x_j \not= z_j\} })}_{1 \leq j \leq n}$, $z \in A$.
If then $\big | \sum \theta_z s'(z) \big|^2$ witnesses $d_A(y)^2$ (up to some $\delta >0$),
$$
d_A(x)^2  - d_A(y)^2 
		\, \leq \, \bigg | \sum \theta_z s(z) \bigg|^2 - \bigg | \sum \theta_z s'(z) \bigg|^2
$$
where $s(z)$ differs from $s'(z)$ only on the $i$-th coordinate which is
$\mathbbm{1}_{\{x_i \not= z_i\} }$ rather than $\mathbbm{1}_{\{y_i \not= z_i\} }$. That is
$$
d_A(x)^2  - d_A(y)^2 
		\, \leq \, \bigg ( \sum \theta_z \mathbbm{1}_{\{x_i \not= z_i\} } \bigg)^2 
		- \bigg (\sum \theta_z \mathbbm{1}_{\{y_i \not= z_i\} } \bigg)^2 \, \leq \, 1
$$
justifying \eqref {eq.lipschitzsquare}.

For $f = e^{-\lambda F^2}$, $\lambda \geq 0$, the reasoning of the
case $f = e^{\lambda F^2}$ may then be repeated together with the new
information \eqref {eq.lipschitzsquare} to get that
whenever $F_i(x_i) \geq F_i(y_i)$ ($\geq 0$), $i = 1, \ldots, n$,
\beqs \begin {split}
\big [-\lambda F_i(x_i)^2 + \lambda F_i(y_i)^2 \big] 
	\big [ e^{-\lambda F_i(x_i)^2} - e^{-\lambda F_i(y_i)^2}\big] 
	 & \, \leq \, \lambda^2 \big [ F_i(x_i)^2 - F_i(y_i)^2) \big]^2  e^{-\lambda F_i(y_i)^2} \\
	  & \, \leq \, \lambda^2 e^\lambda\big [ F_i(x_i)^2 - F_i(y_i)^2) \big]^2  e^{-\lambda F_i(x_i)^2} \\
	  & \, \leq \, 4 \lambda^2 e^\lambda
	  	 \big (a_i(x) + \varepsilon \big)^2 F_i(x_i)^2  e^{-\lambda F_i(x_i)^2} .\\
\end {split} \eeqs
Arguing as before and letting $\varepsilon \to 0$ yields that
$$
\mathrm {Ent}_P (e^{-\lambda F^2}) 
		\, \leq \, 4 \lambda^2 e^\lambda \int_X F^2 e^{-\lambda F^2} dP
$$
for every $ \lambda \geq 0$.

Set $\Lambda (\lambda) = \int_X e^{-\lambda F^2} dP$, $\lambda \in \R$, so that the preceding
expresses that, in the range $0 \leq \lambda \leq \frac 12$ (for example),
$$
  \lambda \Lambda '(\lambda) - \Lambda (\lambda) \log \Lambda (\lambda)
 		\, \leq \, - 8\lambda^2 \Lambda '(\lambda) .
$$
This differential inequality is integrated as before, this time with
$K(\lambda ) = - \log \Lambda (\lambda) \geq 0$, to produce that
$$
- K(\lambda ) = \log \bigg (\int_X e^{-\lambda F^2} dP\bigg)
		 \, \leq \, - \frac {\lambda M_2}{1 + 8\lambda} 
$$
for $0 \leq \lambda \leq \frac 12$.

Since $F = F_A = 0$ on $A$, it follows that for $\lambda = \frac 12$,
$ \frac {1}{10} M_2 \leq    \log \big ( \frac {1}{P(A)} \big)$. Together with
\eqref {eq.convexhullintermediate1},
$$
\int_X e^{ \frac {1}{14} F_A^2} dP  \, \leq \, \frac {1}{P(A)}
$$
which, up to the numerical constant, is the announced Talagrand inequality \eqref {eq.convexhull}.

This thereby concludes the proof of the Talagrand convex hull inequality by the entropy method.
It is not clear (see~\cite {BLM13}) whether the constant $\frac 14$ may be reached by this method
(nor than $\frac 14$ is optimal).

\medskip

\textbf {Some notes and references.} Talagrand's convex distance inequality
was first established in \cite {T91}, following an earlier result on the discrete cube \cite {T88, JS91} -- 
the note \cite {JS91} by W.~Johnson and G.~Schechtman
actually motivated M.~Talagrand towards the general formulation
(observe that $d_A(x) \geq \inf_{y \in \mathrm {Conv}(A)} |x-y|$ on the discrete cube $X = \{0,1\}^n$),
with a (rather short) proof going by induction on the dimension together with geometric arguments.
The original statement \cite {T91,T95} actually involves a stronger family of distances rather
than only the quadratic $d_A^2$. 
The proof presented here via the entropy method was put forward by S.~Boucheron, G.~Lugosi and P.~Massart
in \cite {BLM03} and \cite {BLM13}.

Mass transportation proofs of the convex distance inequality have been considered in
\cite {M96a,M96b,DZ96,D97,S00,DZ10,GRST17}... (see  \cite [Chapter 8] {BLM13} for an account).

Talagrand's convex distance inequality was initially 
motivated by several issues in probability in Banach spaces, in particular \cite {T88} the analogue
of the Gaussian concentration inequality for norms of series
$S = \sum_{i=1}^n \varepsilon_i v_i $ of independent
Bernoulli random variables $\varepsilon_i$,
$\PP ( \varepsilon_i = 1) = \PP (\varepsilon_i = 0) = \frac 12$,
with vector-valued coefficients $v_i$, $i = 1, \ldots , n$,
$$
\PP \Big ( \big |  \| S \| -  \E \big ( \| S \| \big)  \big | \geq r \Big)
	\, \leq \, 4 \, e^{-r^2/4\sigma^2}, \quad r \geq 0,
$$
where $\sigma^2 = \sup _{\|\xi \|\leq 1} \sum_{i=1}^n \langle \xi, v_i \rangle ^2$
(cf.~\cite {LT91}), a significant strengthening of the famous Khintchine-Kahane inequality \cite {K68}.

Its abstract and powerful potential was then emphasized in the monumental
memoir \cite {T95} on concentration inequalities for product measures, describing
numerous illustrations in discrete and combinatorial probability theory. References
to this work give an idea of the impact of the result. It may be refereed in particular
to the general reviews and books \cite {T95,S97,M98,L01,MR02,DP09,T12, V12,BLM13, AS16,V18}...
for a sample applications and developments.

Among the commonly used forms of the convex distance inequality, one may put forward
the following two.

\begin {corollary}
Let $P = \mu_1 \otimes \cdots \otimes \mu_n$ be a product probability measure
on a product space $X = \Omega_1 \times \cdots \times \Omega_n$.
Let $F : X \to \R$ (measurable) be such that for every
$x \in X$ there exists $ a(x) = a  \in [0,\infty)^n$ with $|a| =1$ such
that for every $y \in X$,
\beq \label {eq.lipschitzconvex}
F(x) \, \leq \, F(y) + d_a(x,y). 
\eeq
Then, if $M$ is a median of $F$ for $P$, for any $r \geq 0$,
$$ 
P\big (|F - M| \geq r \big )  \, \leq \,  4\,e^{-r^2/4} .
$$
\end {corollary}

Replacing $F$ by $-F$, this corollary applies similarly if \eqref {eq.lipschitzconvex} is changed into
$ F(y) \leq F(x) + d_a(x,y)$.

\begin {corollary}
Let $P$ be any product probability measure supported
on $[0,1]^n$. For every convex Lipschitz function $F$ on $\R^n$
with ${\|F\|}_{\mathrm {Lip}} \leq 1$, and every $r\geq 0$,
$$ 
P\big (|F - M| \geq r \big )  \, \leq \,  4\,e^{-r^2/4} 
$$
where $M$ is a median of $F$ for $P$.
\end {corollary}

At the expense of numerical constants, medians may be replaced by expectations. Integrating
in $r \geq 0$ the inequalities of the preceding corollaries yields
$| \int_X F dP - M| \leq 4 \sqrt \pi \leq 8$. Hence
$$
	P\big (|F - \textstyle {\int_X F dP}| \geq r \big )  \, \leq \, 
	P\big (|F - M| \geq \textstyle {\frac r2} \big )  \, \leq \,  4\,e^{-r^2/16} 
$$
if $r \geq 16$, while
$ P\big (|F - \int_X F dP| \geq r \big )  \leq 1\leq  e^{16} \, e^{-r^2/16}$
if $r \leq 16$. It is also simple to go back from a concentration inequality around the mean
to one around a median (cf.~\cite {MS86,L01}).

\section {The $\mathrm {L}^1$--$\mathrm {L}^2$ variance inequality} \label {sec.5}

With respect to the original formulation of
Talagrand's $\mathrm {L}^1$--$\mathrm {L}^2$ variance inequality in Section~\ref {sec.1}, 
the argument is developed here
with functions on $\{-1,+1\}^n$ and makes use of the framework presented in Section~\ref {sec.3}.
The end of the proof catches up with the original statement \eqref {eq.l1l2}.
The main point of the proof will be to use
hypercontractivity on the expansion of the variance along the semigroup
${(P_t)}_{t\geq 0}$ defined in \eqref {eq.semigroupcube}. 

The starting point is therefore the variance representation along
${(P_t)}_{t\geq 0}$ of a function $f $ on $X = \{-1,+1\}^n$ as
\beqs \begin {split}  
\int_X f^2 d\mu_{\mathsf{p}}^n - \int_X (P_1 f)^2  d\mu_{\mathsf{p}}^n 
       & \, = \,  - \int_0^1 \bigg (\frac {d}{dt} \int_X (P_t f)^2 d\mu_{\mathsf{p}}^n \bigg ) dt \\
       & \, = \,   - 2 \int_0^1 \bigg ( \int_X P_t f \, \mathrm {L} P_t f d\mu_{\mathsf{p}}^n \bigg ) dt \\
       & \, = \,    2 \int_0^1 \sum_{i=1}^n \int_X  (\mathrm {L}_i P_tf)^2 d\mu_{\mathsf{p}}^n \, dt . \\
\end {split} \eeqs
Assume next that $\int_X f d\mu_{\mathsf{p}}^n = 0$ so that 
$\int_X P_sf d\mu_{\mathsf{p}}^n = 0$ for every $s \geq 0$ as well.
From the Poincaré inequality \eqref {eq.poincare},
the derivative of the map $s \mapsto e^{2s} \int_X (P_s f)^2  d\mu_{\mathsf{p}}^n $
is negative, so the map is decreasing and thus
$$
\int_X (P_1 f)^2  d\mu_{\mathsf{p}}^n  \, \leq \, \frac {1}{e^2} \int_X f^2  d\mu_{\mathsf{p}}^n .
$$
Therefore
$$
\int_X f^2 d\mu_{\mathsf{p}}^n \, \leq \, 
		3 \int_0^1 \sum_{i=1}^n \int_X  (\mathrm {L}_i P_tf)^2 d\mu_{\mathsf{p}}^n \, dt.
$$

Now $\mathrm {L}_i \mathrm {L} f = \mathrm {L}\, \mathrm {L}_i f$ so that
$\mathrm {L}_i P_tf = P_t (\mathrm {L}_i f)$ for every $i=1, \ldots, n$ and $t \geq 0$. 
It may thus be called on the hypercontractivity
property \eqref {eq.hypercontractivitycube} to get that, for every $i = 1, \ldots, n$ and $t>0$,
$$ 
\int_X  (\mathrm {L}_i P_tf)^2 d\mu_{\mathsf{p}}^n 
	\, = \, \int_X  \big | P_t (\mathrm {L}_i f) \big|^2 d\mu_{\mathsf{p}}^n 
		\, \leq \, \bigg (\int_X  |\mathrm {L}_i f|^p  d\mu_{\mathsf{p}}^n\bigg)^{2/p}
$$   
where $ p = p(t) = 1 + e^{-4\rho t} < 2$. Recall that
$ \rho = \frac {\mathsf{p} - \mathsf{q}}{\log \mathsf{p} - \log \mathsf{q}}$
($=\frac 12$ if $\mathsf{p = \frac 12}$).
After the change of variables $ p(t) = v $, it holds that
$$
\int_X f^2 d\mu_{\mathsf{p}}^n
	\, \leq \,  \frac 1\rho \, e^{4\rho}
	 \sum_{i=1}^n  \int_1^2  \bigg (\int_X  |\mathrm {L}_i f|^v d\mu_{\mathsf{p}}^n \bigg)^{2/v} dv. 
$$

This inequality actually basically amounts to the result
(and may be used toward an Orlicz space formulation -- see the comment below). 
Indeed, by H\"older's inequality,
$$ 
  \bigg (\int_X  |\mathrm {L}_i f|^v d\mu_{\mathsf{p}}^n \bigg)^{1/v}  
  	\, = \, {\| \mathrm {L}_i f \|}_v 
 	\, \leq \,    {\| \mathrm {L}_i f \|}_1^\theta  \, 
		   {\| \mathrm {L}_i f \|}_2^{2/\theta}  
$$
where $\theta  = \theta (v) \in [0,1]$ is defined by 
$ \frac 1v = \frac \theta1  + \frac {1-\theta}{2}$.
Hence
$$  
  \int_1^2 {\| \mathrm {L}_i f \|}_v^2 \, dv
    \,  \leq \,   { \|   \mathrm {L}_i f  \|} ^2_2  \int _1^2  b^{2 \theta (v) } dv
$$
where $ b = \frac {{\|  \mathrm {L}_i f\| }_1}{{\|  \mathrm {L}_i f\| }_ 2} \leq 1$
($=1$ if ${\|  \mathrm {L}_i f\| }_ 2 = 0$). It remains to
evaluate the latter integral with $2\theta (v) = s$,
$$ 
\int_1^2 b ^{2 \theta (v) } dv  \, \leq \, \int_0^2 b^s ds \, \leq \,  \frac{2}{1 + \log (\frac 1b)} \, .
$$
As a consequence, for any $f: X \to \R$,
\beq \label {eq.l1l2bis}
	\mathrm {Var}_{\mu_{\mathsf{p}}^n} (f)
	\, \leq \, \frac 2\rho \, e^{4\rho}
	 \sum_{i=1}^n 
	 \frac { { \|   \mathrm {L}_i f  \|} ^2_2}{ 1 
	 + \log \frac {{\|  \mathrm {L}_i f\| }_2}{{\|  \mathrm {L}_i f\| }_ 1}} \, .
\eeq 

It remains to make the comparison with the formulation \eqref {eq.l1l2}
of Talagrand's inequality. To this task, observe that
$$
\mathrm {L}_if(x) 
		\,= \, (1 - \mathsf{p}) 
			\big [ f(x_1, \ldots, x_{i-1}, -1, x_{i+1}, \ldots, x_n) - f(x) \big]
$$
if $x_i = +1$ and
$$
\mathrm {L}_if(x) 
	 \,= \, \mathsf{p} \big [ f(x_1, \ldots, x_{i-1}, +1, x_{i+1}, \ldots, x_n) - f(x) \big]
$$
if $x_i = -1$, that is what is denoted $-\Delta_i f$ in \eqref {eq.l1l2} 
after the change from $\{-1,+1\}$ to $\{0,1\}$. The Talagrand inequality
\eqref {eq.l1l2} therefore follows having observed that $\frac 1\rho \, e^{ 4\rho}$ is of the order
of $\log \frac {1}{\mathsf{p}}$ as $\mathsf{p} \to 0$. In fact, from the preceding proof, $K= 30$
($K=14$ if $\mathsf {p = \frac 12}$) is a valid
numerical constant for \eqref {eq.l1l2} (but may be easily improved).

\medskip

\textbf {Some notes and references.}
The original proof in \cite {T94} uses some Fourier
analysis on the discrete cube (although the author is claiming that it may not be
used in the $\mathsf {p \not= \frac 12}$ case).
It does not mention hypercontractivity although it is implicit
(Lemma~2.1). A simplified proof, for $\mathsf {p = \frac 12}$, was proposed
in \cite {BKS03} putting forward the hypercontractive argument. 
The proof presented here is taken from \cite {CEL12}, where it is
extended to wider settings. In particular, the Talagrand $\mathrm {L}^1$--$\mathrm {L}^2$
inequality also applies to the Gaussian model.

Talagrand's $\mathrm {L}^1$--$\mathrm {L}^2$ variance
inequality was motivated by a result of L.~Russo \cite {R82} on a threshold effect for monotone sets
depending little on any given coordinate. The result also provided an alternate proof of the
famous result of J.~Kahn, G.~Kalai
and N.~Linial \cite {KKL88} about influences on the cube (that already used
hypercontractivity). Actually, M.~Talagrand's approach is an adaptation of the ideas of \cite {KKL88}.
Namely, applying \eqref {eq.l1l2} to the (Boolean)
function $f = \mathbbm {1}_A - \mu(A) $ for some set $A \subset \{0, 1\}^n$ with
$\mu_{\mathsf{\frac 12}}^n (A) = \alpha \in (0,1)$ (take $\mathsf {p=\frac 12}$ for simplicity), it follows that
$$ 
\alpha (1-\alpha) \,  \leq  \, 2K \sum_{i=1}^n
		 \frac {\mathrm {I}_i (A)}{ \log \big (\frac {e}{ \sqrt {2\mathrm {I}_i (A)}} \, \big )} 
$$
where, for each $i = 1, \ldots, n$,
$$ 
 \mathrm {I}_i (A) \, = \,  \mu_{\mathsf{\frac 12}}^n \big ( x \in A ;  U_i (x) \notin A  \big )
$$
is the so-called influence of the $i$-th coordinate on the set $A$. 
In particular, there is a coordinate $i$, $1 \leq i\leq n$, such that
$$ 
	\mathrm {I}_i (A)   \, \geq \,  \frac{\alpha(1-\alpha)}{8K} \, \frac {\log n}{n}
$$
which is the main result of \cite {KKL88}. This result
remarkably improves by a (optimal) factor $\log n$ what would
follow from the Poincaré inequality \eqref {eq.poincare} applied to $ f = \mathbbm {1}_A$.
Since then, the Talagrand $\mathrm {L}^1$--$\mathrm {L}^2$ variance inequality
plays a major role in Boolean analysis and its ramifications with theoretical
computer science (cf.~e.g.~\cite {OD14}).

The Talagrand inequality \eqref {eq.l1l2} indeed
represents a sharpening upon the Poincaré inequality \eqref {eq.poincare}
(up to numerical constants).
It may actually be interpreted as a dual form of the logarithmic Sobolev inequality
in the (loose) sense that the latter ensures that if some gradient of a function $f$ is 
in $\mathrm {L}^2$, then the function belongs to the Orlicz space $\mathrm {L}^2 \log \mathrm {L}$,
while the Talagrand inequality expresses that if the gradient is in
$\mathrm {L}^2 (\log \mathrm {L})^{-1}$, then the function (with zero mean) is in $\mathrm {L}^2$
(cf.~Theorem 1.6 in \cite {T94}).
This dual point of view was emphasized in \cite {BH99}. A specific feature of the 
Talagrand inequality \eqref {eq.l1l2} is nevertheless that is has a suitable product structure.
Alternate forms of the Talagrand inequality as a direct consequence of the
logarithmic Sobolev inequality are developed in \cite {FS07} and \cite {R06}.

With the work \cite {BKS03} by I.~Benjamini, G.~Kalai and O.~Schramm (see further [BR08]),
the Talagrand $\mathrm {L}^1$--$\mathrm {L}^2$ variance
inequality has been identified as one of the rare tool towards sub-diffusive
regimes and super-concentration phenomena, ubiquitous to many models of the current research
(percolation, random matrices, spin glasses etc.). The surveys and monographs
\cite {BKS99,C14,GS15,ADH17,S18}... give an account on these recent active developments.

\section {The quadratic transportation cost inequality} \label {sec.6}

This section presents the entropic proof of Talagrand's transportation inequality \eqref {eq.transport},
deducing it from the logarithmic Sobolev inequality \eqref {eq.logsobgaussian}.

To start with, recast the transportation metric $T_w$ of \eqref {eq.transport} as 
the (quadratic) Kantorovich distance between probability measures $\mu$ and $\nu$ 
on the Borel sets of $\R^n$
$$
\mathrm {W}_2 (\mu, \nu)
		 \, = \, \inf \bigg ( \int_{\R^n \times \R^n} |x-y|^2 d\pi (x,y) \bigg)^{1/2}
$$
where the infimum is taken over all couplings $\pi$ on $\R^n \times \R^n$ with respective
marginals $\mu$ and $\nu$, $|x-y|$ being the Euclidean distance between $x$ and $y$ in $\R^n$.

The classical duality formula (e.g.~\cite {V03}) for the Kantorovich distance $\mathrm {W}_2(\mu, \nu)$
between two probability measures $\mu$ and $\nu$ on the Borel sets of $\R^n$ expresses that
\beq \label {eq.w2dual}
 \frac {1}{2} \, \mathrm{W}_2(\mu ,\nu ) ^2
    \, = \, \sup \bigg ( \int_{\R^n}  Q_1 \varphi  \,  d\mu  -  \int _{\R^n} \varphi  \,  d\nu  \bigg )
\eeq
where the supremum is taken over all bounded continuous functions $\varphi  : \R^n \to \R$
and where
$$
 Q_s \varphi  (x) \, = \,  \inf_{y \in \R^n} \Big [ \varphi (y) + \frac{1}{2s} \, |x-y|^2 \Big],
   \quad  s>0, \, \, x \in \R^n, 
$$
is the infimum-convolution Hopf-Lax semigroup. It is standard (cf.~e.g.~\cite{E98,V03}) that
$Q_s \varphi (x)$, $s>0$, $x \in \R^n$, solves the Halmiton-Jacobi equation
\beq \label {eq.hamiltonjacobi}
\frac {d}{ds} \, Q_s \varphi \, = \, -\frac 12 \, |\nabla Q_s \varphi |^2 
\eeq
in $(0,\infty) \times \R^n$ with initial condition $\varphi$.

On the basis of \eqref {eq.w2dual}, the
entropic inequality \eqref {eq.entropicinequality} ensures that if $d\mu = f d\nu$,
\beqs \begin {split}
\int_{\R^n}  Q_1 \varphi  \,  d\mu  -  \int _{\R^n} \varphi  \,  d\nu 
	& \, = \, \int_{\R^n}  Q_1 \varphi  \, f \, d\nu  -  \int _{\R^n} \varphi  \,  d\nu  \\
	& \, \leq \, \int_{\R^n} f \log f \, d \nu + \log \int_{\R^n} e^{Q_1 \varphi} d\nu
				-  \int _{\R^n} \varphi  \,  d\nu . \\
\end {split} \eeqs

Assume now that $\nu$ is the standard Gaussian measure $\gamma$ on $\R^n$. 
The proof will follow the Herbst argument and deduce the Talagrand transportation cost
inequality from the logarithmic Sobolev inequality \eqref {eq.logsobgaussian}. Apply namely the latter
$f^2 =  e^{s \, Q_s \varphi}$, $s >0$, to get that
$$
 \int_{\R^n} s \, Q_s \varphi \, e^{s \, Q_s \varphi} d\gamma - \Lambda(s) \log \Lambda(s)
		\, \leq \, \frac {s^2}{2} \int_{\R^n} |\nabla Q_s \varphi |^2 e^{s \, Q_s \varphi} d\gamma
$$
where $ \Lambda (s)= \int_{\R^n} e^{s \, Q_s \varphi} d\gamma$.
But $ \partial_s Q_s \varphi = -\frac 12 |\nabla Q_s \varphi |^2 $ so that
$$
\Lambda'(s) \, = \,  \int_{\R^n} Q_s \varphi \, e^{s \, Q_s \varphi} d\gamma
		- \frac s2  \int_{\R^n} |\nabla Q_s \varphi |^2 e^{s \, Q_s \varphi} d\gamma,
$$
and hence the previous inequality is translated into
$$
s \, \Lambda' (s) - \Lambda (s) \log \Lambda (s) \, \leq \, 0
$$
for every $s>0$.
Setting $H(s) = \frac 1s \, \log \Lambda(s)$, $H(0) = \int_{\R^n} \varphi \, d\gamma$, it 
therefore holds true that $H'(s) \leq 0$, $s >0$. Hence
$$
 \log \int_{\R^n} e^{Q_1 \varphi} d\gamma \, = \, H(1) \, \leq \, H(0) \, = \, \int_{\R^n} \varphi \, d\gamma ,
$$
and plugging this conclusion into the above entropic inequality yields that
$$
\int_{\R^n}  Q_1 \varphi  \,  d\mu  -  \int _{\R^n} \varphi  \,  d\gamma 
		\, \leq \, \int_{\R^n} f \log f \, d \gamma.
$$
Taking the supremum over all bounded continuous $\varphi : \R^n \to \R$,
the Kantorovich duality \eqref {eq.w2dual} yields
$$
\mathrm{W}_2 (\mu ,\gamma )^2 \, \leq \, 2\int_{\R^n} f \log f \, d \gamma
$$
which is the Talagrand quadratic transportation cost inequality \eqref {eq.transport}.

As for the logarithmic Sobolev inequality \eqref {eq.logsobgaussian}, the constant $2$
is optimal in \eqref {eq.transport} as can be seen by the choice for $\mu$ of a shift of $\gamma$
by $b \in \R^n$.

\medskip

\textbf {Some notes and references.} M. Talagrand's original proof of \eqref {eq.transport}
in \cite {T96a} uses monotone transport in dimension one together with a tensorization argument.
The main result of \cite {T96a} is actually a corresponding
(stronger) inequality for products of the exponential measure, the Gaussian case being presented
as a simpler example to deal with first. Motivation comes from the mass transportation approach
to the concentration of measure phenomenon put forward by K.~Marton \cite {M86,M96a},
and a sharp form of which for these measures.

Mass transportation proofs directly in dimension $n$ have been 
provided next, in particular by means of the Brenier map \cite {CE02,B03,V03,V09}...
That the logarithmic Sobolev inequality implies the Talagrand quadratic transportation cost
inequality is the main achievement of the celebrated paper \cite {OV00} by F.~Otto and
C.~Villani. The implication holds in a rather general setting. The approach in \cite {OV00}
relies on the formal Otto calculus in Wasserstein space \cite {O01}, and has been one
driving force in the study of functional inequalities and curvature lower bounds in
metric measure spaces (cf.~\cite {V03,V09, BGL14}). The proof presented here inspired
by the Herbst argument is taken from \cite {BGL01,BG99}. For more on transportation cost
inequalities, see in addition \cite {GL10}.

That mass transportation is at the root of the four Talagrand inequalities is an
alternate project of independent interest, already partly understaken in \cite [Chapter 8]{BLM13}.

\section {Supremum of empirical processes} \label {sec.7}

This last section is thus devoted to the proof of the inequality \eqref {eq.supremum},
in the following notation according to the choices
adopted so far. Let $X_1, \ldots, X_n$ be independent random variables
on a probability space $(\Omega, \mathcal {A}, \PP)$ with values in some
measurable space $(S, \mathcal {S})$. Let $\mathcal {F}$ be a countable family of
real-valued uniformly bounded measurable functions on $S$, and set
\beq \label {eq.z}
Z \, = \, \sup_{g \in \mathcal {F}} \sum_{i=1}^n g(X_i).
\eeq
By (dominated) convergence and homogeneity, it is enough to consider a finite family
$\mathcal {F} = \{g_1, \ldots, g_N\}$ such that $|g_k| \leq 1$, $k = 1, \ldots, N$.
(It is also assumed below that $Z$ is not $0$ almost surely.)

To make use of the framework of Section~\ref {sec.3}, denote
by $\mu_1, \dots, \mu_n$ the respective probability distributions of the
independent random variables $X_1, \ldots, X_n$,
and set $P = \mu_1 \otimes \cdots \otimes \mu_n$ on the product space $S^n$. In accordance
$$
Z \, = \, Z(x) \, = \, Z(x_1, \ldots, x_n) \, = \, \max_{1 \leq k \leq N} \sum_{i=1}^n g_k(x_i),
	\quad x = (x_1, \ldots, x_n) \in S^n
$$
(for which nevertheless the probabilistic notation induced by
\eqref {eq.z} will be used from time to time below). Following the Herbst argument, 
the task will be to apply the tensorized logarithmic Sobolev inequalities of
\eqref {eq.tensorization1} and
\eqref {eq.tensorization2} to $e^{\lambda Z}$ for every $\lambda \in \R$,
in much the same way as for the convex distance inequality in Section~\ref {sec.4}.

As is classical in the study of exponential inequalities for sums of independent random
variables, the inequality \eqref {eq.supremum} entails a Gaussian tail for the small
values of $r$ and a Poisson one for the large values.

In the first part of the proof, consider Poisson tails for non-negative random variables.
Let thus $g_k$, $k = 1, \ldots, N$, be such
that $0 \leq g_k \leq 1$ and apply \eqref {eq.tensorization2} to $f = e^{\lambda Z}$,
$\lambda \in \R$. For each $i = 1, \ldots, n$ and $x = (x_1, \ldots, x_n) \in S^n$, choose
$c_i = e^{\lambda Z^i(x)}$ where
$$
Z^i (x) \, = \,  \max_{1 \leq k \leq N} \sum_{j \not= i} g_k(x_j).
$$
By definition $Z^i (x)$ only depends on $x_1, \ldots, x_{i-1}, x_{i+1}, \ldots, x_n$ and
$0 \leq Z_i(x_i) - Z^i(x) \leq  1$, $i=1, \ldots, n$.
Now
$$
	\lambda \big [Z_i(x_i)- Z^i (x) \big] e^{\lambda Z_i(x_i)}  - [ e^{\lambda Z_i(x_i)} - e^{\lambda Z^i(x)} ]
	 \, = \, \phi \big ( - \lambda [ Z_i(x_i) - Z^i (x)] \big) \, e^{\lambda Z_i(x_i)}
$$
where $\phi (u) = e^u - 1 - u$, $u \in \R$. Since $\phi$ is convex and $\phi (0)=0$,
$\phi (- \lambda u ) \leq u \phi (-\lambda)$ for every $\lambda $ and $0 \leq u \leq 1$, so that
\beqs \begin {split}
  \sum_{i=1}^n
	\lambda \big [Z_i(x_i)- Z^i (x) \big] e^{\lambda Z_i(x_i)} 
		& - [ e^{\lambda Z_i(x_i)} - e^{\lambda Z^i(x)} ]  \\
	& \, \leq \,  \sum_{i=1}^n \big [ Z_i(x_i) - Z^i (x) \big] \phi (-\lambda ) \, e^{\lambda Z_i(x_i)}\\
	& \, \leq \, \phi (-\lambda ) \, Z(x) \, e^{\lambda Z(x)}  \\
\end {split} \eeqs
where it is used that
$$
\sum_{i=1}^n \big [ Z_i(x_i) - Z^i (x) \big] \, = \, 
		\sum_{i=1}^n  \big [ Z(x) - Z^i (x) \big ]  \, \leq \, Z(x).
$$

As a consequence therefore of \eqref {eq.tensorization2}, and
with probabilistic notation, for any $\lambda \in \R$,
$$
	\lambda \, \E ( Z e^{\lambda Z} ) -  \E ( e^{\lambda Z} )   \log \E ( e^{\lambda Z} )
	 \, \leq \,  \phi (-\lambda) \,  \E ( Z e^{\lambda Z} ).
$$ 
If $\Lambda (\lambda) = \E ( e^{\lambda Z} )$ and 
$H(\lambda) =  \frac 1\lambda \log \Lambda(\lambda)$, $\lambda \in \R$
($H(0) = \int_X Z dP$), the preceding reads
$$
H'(\lambda ) \, \leq \, \frac {\phi (-\lambda)}{\lambda^2} 
		\, \frac {\Lambda' (\lambda)}{\Lambda (\lambda)} \, .
$$
Since $H'(\lambda) 
= - \frac 1\lambda H(\lambda) + \frac 1 \lambda \frac {\Lambda' (\lambda)}{\Lambda (\lambda)} $,
it follows that for $\lambda >0$,
\beq \label {eq.herbstphi}
 \frac {H'(\lambda)}{H(\lambda)} 
 	\, \leq \,  \frac {\phi (-\lambda)}{\lambda (\lambda - \phi (-\lambda))} 
	\, = \,  \frac {1}{\lambda - \phi (-\lambda)} - \frac 1\lambda \, .
\eeq
It is easily seen that, again with $\lambda >0$,
$$
\int_0^\lambda  \Big [ \frac {1}{u - \phi (- u)}  - \frac 1u \Big ] \, du
	\, = \, \int_0^\lambda \Big [ \frac {1}{1 - e^{-u}}  - \frac 1u \Big ] \, du
	\, = \, \lambda  -\log \lambda + \log (1 - e^{-\lambda})
$$
so that, integrating \eqref {eq.herbstphi},
$$
H(\lambda ) \, \leq \, H(0) \frac 1 \lambda \, (e^\lambda - 1) , \quad \lambda \geq 0.
$$

As a conclusion of this analysis, the following statement holds true.

\begin {proposition} \label {prop.poisson}
If $0 \leq g \leq 1$ for every $g \in \mathcal {F}$,
$$
\E(e^{\lambda Z} ) \, \leq \, e^{\E(Z)(e^\lambda - 1)}
$$
for every $\lambda \geq 0$. As a consequence, for every $r\geq 0$,
\beq \label {eq.poissontail}
\PP \big ( Z \geq \E(Z) + r \big)
	\, \leq \, \exp \bigg ( - \E (Z) \, h \Big ( \frac {r}{\E(Z)} \Big) \bigg)
\eeq
where $h(u) = (1+u) \log (1+u) - u$, $u \geq 0$.
\end {proposition}

The Poisson tail \eqref {eq.poissontail} is obtained from Markov's inequality
$$
\PP \big ( Z \geq \E(Z) + r \big) \, \leq \, 
	e^{-\lambda (\E(Z) + r) + \E(Z)(e^\lambda - 1)} 
$$
and optimization in $\lambda \geq 0$.

For functions $g$ taking values in $[-1,+1]$,
a similar scheme may be followed on the basis this time of \eqref {eq.tensorization1} towards Gaussian tails.
Let thus $g_k$, $k = 1, \ldots, N$, be such
that $ |g_k| \leq 1$ and apply \eqref {eq.tensorization1} to $f = e^{\lambda Z}$.
If $Z_i(x_i) \geq Z_i(y_i)$, $x_i, y_i \in S$, $i=1, \ldots, n$, and $\lambda \geq 0$,
by the mean value theorem,
$$
	 \lambda \big [ Z_i(x_i) - Z_i(y_i)\big]
	 \big [ e^{\lambda Z_i(x_i)} - e^{\lambda Z_i(y_i)} \big]
	  \, \leq \, \lambda^2 \big [ Z_i(x_i) - Z_i(y_i)\big]^2 e^{\lambda Z_i(x_i)} .
$$	
By the definition of $Z=Z(x)$,
$$
\sum_{i=1}^n \big [ Z_i(x_i) -  Z_i (y_i)\big ]^2 \mathbbm {1}_{\{ Z_i(x_i) \geq Z_i(y_i) \}}
		\, \leq \, \max_{1 \leq k \leq N}  \sum_{i=1}^n \big [ g_k(x_i) - g_k (y_i) \big]^2.
$$
Denoting by $\widetilde {W} = \widetilde {W}(x,y)$ the right-hand side of the preceding inequality, 
\eqref {eq.tensorization1} indicates that, in probabilistic notation,
\beq \label {eq.empirical+}
 \lambda \, \E (Z e^{\lambda Z} ) - \E ( e^{\lambda Z} ) \log \E ( e^{\lambda Z} )
		\, \leq \, \lambda^2 \, \E ( \widetilde {W} e^{\lambda Z})
\eeq
for any $\lambda \geq 0$.

Working with $\lambda \leq 0$, it still holds that
$$
	 \lambda \big [ Z_i(x_i) - Z_i(y_i)\big]
	 \big [ e^{\lambda Z_i(x_i)} - e^{\lambda Z_i(y_i)} \big]
	  \, \leq \, \lambda^2 e^{-2\lambda}\big [ Z_i(x_i) - Z_i(y_i)\big]^2 e^{\lambda Z_i(x_i)}
$$
as $ (0 \leq) \, Z_i(x_i) -  Z_i (y_i) \leq 2$. Together with \eqref {eq.empirical+}, it may thus be
concluded that, in the range $\lambda \in [-\frac 14,+\frac 14]$ (for example),
$$
 \lambda \, \E (Z e^{\lambda Z} ) - \E ( e^{\lambda Z} ) \log \E ( e^{\lambda Z} )
		\, \leq \, 2 \lambda^2 \, \E ( \widetilde {W} e^{\lambda Z}).
$$
In addition, by independence, 
$\E ( \widetilde {W} e^{\lambda Z}) \leq 2V\, \E(e^{\lambda Z}) + 2\, \E ( W e^{\lambda Z})$
where $ W = W(x) = \max_{1 \leq k \leq N}  \sum_{i=1}^n g_k(x_i)^2$ and, in the notation
of \eqref {eq.supremum}, $V = \E(W)$. As a consequence, for every
$\lambda \in [-\frac 14,+\frac 14]$, 
\beq \label {eq.empirical+-}
 \lambda \, \E (Z e^{\lambda Z} ) - \E ( e^{\lambda Z} ) \log \E ( e^{\lambda Z} )
		\, \leq \, 4 V \lambda^2 \, \E ( e^{\lambda Z} )+ 4\lambda^2 \,\E ( W e^{\lambda Z}).
\eeq

As in the Herbst argument, the latter inequality \eqref {eq.empirical+-} is transformed
into a differential inequality on Laplace transforms. To start with, observe
that for any $\lambda \in \R$,
\beq \label {eq.entropicw}
\lambda \, \E ( W e^{\lambda Z})
	\, \leq \, \lambda \, \E ( Z e^{\lambda Z} )
	 -  \E (e^{\lambda Z})   \log \E ( e^{\lambda Z}) + \E( e^{\lambda Z} ) \log \E( e^{\lambda W} )
\eeq
(as another instance of the entropic inequality \eqref {eq.entropicinequality}
with $f = \frac {e^{\lambda Z}}{\E(e^{\lambda Z})}$ and $g = \lambda W$).

Start next with the positive values of $\lambda$.
Setting $\Lambda (\lambda) = \E( e^{\lambda Z} ) $, $R (\lambda) = \E(   e^{\lambda W} ) $,
$\lambda \in \R$, it follows from \eqref {eq.empirical+-} and \eqref {eq.entropicw} that for every
$\lambda \in [0,\frac 14]$,
$$
	(1- 4\lambda)
	\big [ \lambda \Lambda'(\lambda) - \Lambda (\lambda) \log \Lambda (\lambda) \big]
		\, \leq \,  4V\lambda^2 \Lambda(\lambda)
			+ 4\lambda \Lambda (\lambda) \log R(\lambda).
$$
Proposition~\ref {prop.poisson} may be applied to $ W$ to get that
$$
\log R(\lambda) \, \leq \, \E(W) (e^{\lambda} - 1)
		\, = \, V (e^{\lambda} - 1)  , \quad \lambda \geq 0.
$$
Hence
$$
(1 - 4\lambda) \big [ \lambda \Lambda'(\lambda) - \Lambda (\lambda) \log \Lambda (\lambda) \big]
	\, \leq \,  4V \big [ \lambda^2 + \lambda (e^\lambda - 1) \big] \Lambda (\lambda) ,
$$
and in the usual notation $H(\lambda) = \frac 1\lambda \log \Lambda (\lambda)$,
for any $0 \leq \lambda < \frac 14$,
$$
H'(\lambda) \, \leq \, \frac {4V}{1 - 4\lambda} 
	\bigg [1 + \frac {e^{\lambda} -1}{\lambda } \bigg].
$$
It remains to suitably integrate this differential inequality. Without trying sharp bounds,
the right-hand side may simply be upper-bounded by $20V$ in the range
$0 \leq \lambda \leq \frac {1}{8}$ so to get that
$$
\E(e^{\lambda Z}) \, \leq \, e^{\lambda \E(Z) + 20\lambda^2 V}.
$$

Given then $r \geq 0$, use Markov's exponential inequality with
$\lambda = \frac {r}{40V}$ if $ r \leq 5V$ and $\lambda = \frac {1}{8}$
otherwise to derive that
$$
\PP \big ( Z \geq \E(Z) + r ) \, \leq \, e^{- \min ( \frac {r}{16} , \frac{ r^2}{80V} )} .
$$

Now, the same reasoning may be applied to $-Z$ since \eqref {eq.empirical+-} holds also for 
$ \lambda \in [-\frac {1}{4}, 0]$. Together with the two parts and the union bound,
the following statement follows.

\begin {proposition} \label {prop.gaussian}
In the preceding notation, for every $r \geq 0$
$$
\PP \big ( \big | Z - \E(Z) \big| \geq r \big) 
		\, \leq \, 2 \, e^{- \min ( \frac {r}{16} , \frac{ r^2}{80 V})}.
$$
\end {proposition}

This proposition is close to the Talagrand inequality \eqref {eq.supremum} with the Gaussian
tail for the small values of $r$ ($\leq 5V$), but only
an exponential decay for the large values (describing the so-called Bernstein inequality,
cf.~\cite {BLM13}). To reach the Poisson tail
and the full conclusion, it should be combined with Proposition~\ref {prop.poisson}.

To this task, let $\tau >0$ and set
$$
Z_\tau^1 \, = \, \max_{1 \leq k \leq N} \sum_{i=1}^n g_k(X_i) 
		\mathbbm {1}_{\{ | g_k(X_i) | \leq \tau \}}
$$
and
$$
Z_\tau^2 \, = \, \max_{1 \leq k \leq N} \sum_{i=1}^n \big | g_k(X_i) \big |
		\mathbbm {1}_{\{ | g_k(X_i) | > \tau \}} 
$$
so that $|Z - Z_\tau^1| \leq Z_\tau^2$. Then, for $r \geq \E(Z_\tau^2)$,
\beq \label {eq.z1z2}
\PP \big ( \big | Z - \E(Z) \big| \geq 4r \big) 
	\, \leq \, \PP \big ( \big | Z_\tau^1 - \E(Z_\tau^1) \big| \geq r \big) 
			+ \PP \big (  Z_\tau^2 \geq  \E(Z_\tau^2) + r \big) .
\eeq
By Proposition \ref {prop.gaussian} applied to $\frac 1\tau Z_\tau^1$, 
$$
\PP \big ( \big | Z^1_\tau - \E(Z^1_\tau) \big| \geq r \big) 
		\, \leq \, 2 \, e^{- \min ( \frac {r}{16\tau} , \frac{ r^2}{80 \tau^2 V})}
$$
for every $r \geq 0$, while by Proposition \ref {prop.poisson} applied to $Z_\tau^2$,
$$
\PP \big (  Z_\tau^2 \geq \E(Z_\tau^2) + r \big)  
	\, \leq \,  e^{- \frac r2 \log ( 1 + \frac {r}{\E(Z_\tau^2) })}
$$
since $h(u) \geq \frac u2 \, \log (1+u)$, $u \geq 0$.

Choose next $ \tau = \sqrt { \frac {4V}{5r} }$ ($ r>0$). In the range $4r \geq 5 V$,
$$
\E(Z_\tau^2 ) \, \leq \, \frac V\tau  \, = \,  \sqrt {\frac 54 \, rV} \, \leq \,  r.
$$
Hence
$$
\frac r2 \log \bigg( 1 + \frac {r}{\E(Z_\tau^2) } \bigg)
	\, \geq \,  \frac r2 \log \bigg ( 1 + \sqrt {\frac {4r}{5V} } \, \bigg)
	\, \geq \, \frac {r}{6} \log \bigg( 1 + \frac {4r}{V } \bigg) .
$$
It also holds in this range that
$$
\min \bigg ( \frac {r}{16\tau} , \frac{ r^2}{80 \tau^2 V} \bigg)
	\, \geq \, \frac {r}{75} \log \bigg ( 1 + \frac {4r}{V } \bigg ) .
$$
Since by the choice of $\tau$, $r \geq \E(Z_\tau^2)$ whenever
$4r \geq 5V$, as a consequence of \eqref  {eq.z1z2},
$$
\PP \big ( \big | Z - \E(Z) \big| \geq 4r \big) 
	\, \leq \, 3 \, e^{- \frac {r}{75} \log ( 1 + \frac {4r}{V })}.
$$

Proposition \ref {prop.gaussian} for the values of $r\leq 5V$ shows that
$$
\PP \big ( \big | Z - \E(Z) \big| \geq r \big)  
	 \, \leq \, 2 \, e^{- \frac {r^2}{80V}}  \, \leq \, 
	 2 \, e^{- \frac {r}{80} \log ( 1 + \frac {r}{V })}
$$
since $\log (1+u) \leq u$, $u \geq 0$. Together with the previous inequality (after
the change from $4r$ to~$r$) yields finally that
$$
\PP \big ( \big | Z - \E(Z) \big| \geq r \big)  
		\, \leq \, 3 \, \exp \bigg ( - \frac {r}{300} \log \bigg ( 1 + \frac {r}{V } \bigg) \bigg)
$$
for every $r \geq 0$. This is the Talagrand inequality \eqref {eq.supremum}, with
$U=1$ by homogeneity, completing thereby its proof.

\medskip

\textbf {Some notes and references.} The proof of \eqref {eq.supremum} developed
in \cite {T96b} is rather cumbersome, elaborating on and deepening the investigation
\cite {T95}. These strengthenings have been clarified since then
by means of information-theoretic and infimum-convolution tools in
several contributions, including
\cite {M96a,M96b,DZ96,D97,S00,P01,Sa03,S07,DZ10,BLM13,GRST17}...
It is deduced from the convex distance inequality after a symmetrization argument in \cite {P03}.
The approach relying on the entropy method and logarithmic Sobolev
inequality was initiated in \cite {L96b}, and expanded and
made precise in \cite {M00a} and in various subsequent publications (cf.~\cite {BLM13}).
The truncation argument to suitably combine the Gaussian and Poisson tails is already
present in \cite {T96b}.
The article \cite {M00a} by P.~Massart provides a much more careful analysis of the involved
numerical constants, of significant relevance in the applications
(Proposition~\ref {prop.poisson} is taken from there -- it is remarkable that it
is already optimal for a class of functions reduced to one element).
The reference \cite {BLM13} provides an account on these developments and on the various steps
and contributions towards
sharper (sometimes optimal) constants in families of inequalities for empirical processes.
The Talagrand inequality \eqref {eq.supremum} on the supremum of empirical processes is indeed
nowadays a major tool in non-asymptotic statistics, where sensitive numerical constants are of importance.
The monographs \cite {M00b,M07,BLM13,GN15,W19}... and the references
therein present numerous illustrations and applications of this most powerful result
in modern statistics.

It is worthwhile mentioning that for the applications, standard symmetrization tools allow
for the bound
$$
V  \, =  \, \E \bigg ( \sup_{f \in  {\cal F}} \sum _{i=1}^n f(X_i)^2 \bigg )
		\, \leq \, U \, \E (Z) + 8 \, \sup_{f \in \mathcal {F}} \sum_{i=1}^n \E \big ( f(X_i)^2 \big)
$$
whenever $\E(f(X_i)) = 0$, $i=1, \ldots, n$, $ f \in \mathcal {F}$ and the class
$\mathcal {F}$ is symmetric ($- \mathcal {F} = \mathcal {F}$) (see \cite {LT91,BLM13}).

\vskip 8mm

\font\tenrm =cmr10  {\tenrm

\parskip 0mm

\noindent Institut de Math\'ematiques de Toulouse 

\noindent Universit\'e de Toulouse -- Paul-Sabatier, F-31062 Toulouse, France

\noindent \&  Institut Universitaire de France 

\noindent ledoux@math.univ-toulouse.fr

}

\end {document}